\documentclass[11pt, fleqn]{amsart}
\usepackage{amsmath, amssymb, amsthm} 
\usepackage{hyperref}
\usepackage{mathtools}
\usepackage{graphicx}
\usepackage{adjustbox}
\usepackage{subcaption}
\newtheorem{theorem}{Theorem}[section]
\newtheorem{lemma}[theorem]{Lemma}
\newtheorem{corollary}[theorem]{Corollary}
\newtheorem{conjecture}{Conjecture}
\newtheorem{proposition}{Proposition}
\newtheorem*{proposition*}{Proposition}

\theoremstyle{remark}
\newtheorem{remark}{Remark}
\title{Matrix methods for arithmetic functions}
\author{Barry Brent}
\email{axcjh@bu.edu}
\thanks{I am grateful to Robin Chapman, Tushar Das, Jos\'e ~L\'opez-Bonilla, and Aleksandar Petojevi\'{c} for generous advice.}
\begin{document}
\begin{abstract} We apply matrix methods to
the study of functions, such as Ramanujan's function $\tau(n)$, the sum of divisors function $\sigma(n)$, and the partion function $p(n)$, that satisfy a certain additive convolution relation. Much of this article has roots in 
the study of symmetric functions, but we do not rely on that theory here. Typically, given two such functions $f$ and $g$, we find determinant
expressions for $f$ and $g$ by analyzing differential equations satisfied by their generating functions. We derive partition-theoretic identities from additive convolutions relating $f$ and $g$, and we study the characteristic polynomials attached to the determinants associated to
$f$ and $g$ (with help from J. L\'opez-Bonilla.)  On entirely empirical grounds, we indicate the possibility of geometric structure among the roots of 
the characteristic polynomials in question.

\end{abstract}
\maketitle
\vspace{-0.5em}
\begin{center}
    \rm{20 October 2025 no2}
\end{center}
\vspace{1em}
\section{Introduction}
\subsection{Contents of this article}
We apply matrix methods to
the study of pairs of functions, such as Ramanujan's function $\tau(n)$ and the sum of divisors function $\sigma(n)$, or the partition function $p(n)$ and the same sum of divisors function (in the opposite order, this time), that satisfy a certain additive convolution relation. As we will see, given an arbitrary function on the positive integers, the companion function exists and can be computed from the convolution relation we have mentioned.   An evaluation of $\tau(n+1)$  as  $(-1)^n|J_n(\overline{S})|/n!$ for an $n$ by $n$ matrix $J_n(\overline{S})$ built from values of the sum of divisors function $\sigma(k), k = 1, 2, ..., n$, is given. 
The characteristic polynomial of $J_n(\overline{S})$ is described in a conjecture of the  writer  proved 
by J. L\'opez-Bonilla. We state a convolution identity
\(
n p(n) = \sum_{k=1}^n \frac{(-1)^k}{24}|H_k(\overline{T})| p(n-k)
\) in which the matrix \(H_k(\overline{T})\) is 
populated with constants and values of Ramanujan's tau function.
The numbers $\tau(n+1)$ and \(p(n)\) are expressed as sums over integer partitions of monomials in the \(\sigma(k)\),
and  \(p(n)\) is expressed as a sum over partitions of monomials in values of tau. In the final section,
we study the possibility of significant geometric structure among the root sets of the relevant characteristic polynomials, in particular, characteristic polynomials of matrices that generate the values of \(\tau(p_n)\) where \(p_n\) denotes the
\(n^{th}\) prime number.

We collect references to relevant earlier work in the bibliography.
\section{Lemmas}
The first lemma appears in (we believe) an erroneous form on
 page 119 of
Adams'  reference work ~\cite{A22}.  Entry 0.313 of Gradshteyn's and Ryzhik's reference work ~\cite{GR88} propagates Adams' version.
Neither book provides a proof
or a citation to one. 
Both give the denominator of the fraction 
\(c_n\)
as \(a_0^n\), instead of what we think is the 
correct \(a_0^{n+1}\). (We note that
there is another formulation of the coefficient $c_n$ in terms of a different determinant in Gholami's article ~\cite{Gh11}.)
\begin{lemma} [corrected version of Adams \cite{A22}, entry 6.360]
Let $\sum_{k=0}^{\infty} a_k x^k$ and $\sum_{k=0}^{\infty} b_k x^k$ be power series with $a_0 \neq 0$. 
If 
\[
\left(\sum_{k=0}^{\infty} b_k x^k\right)\Bigg/\left(\sum_{k=0}^{\infty} a_k x^k\right)
=
\sum_{k=0}^{\infty} c_k x^k,
\]
 then
$c_n = (-1)^n \left |M_n \right |/a_0^{n+1}$,
where $M_n$ is the $n \times n$ matrix
\begin{align}
M_n = \begin{pmatrix}
a_1 b_0 - a_0 b_1 & a_0 & 0 & \cdots & 0 \\
a_2 b_0 - a_0 b_2 & a_1 & a_0 & \cdots & 0 \\
a_3 b_0 - a_0 b_3 & a_2 & a_1 & \cdots & 0 \\
\vdots & \vdots & \vdots & \ddots & \vdots \\
a_n b_0 - a_0 b_n & a_{n-1} & a_{n-2} & \cdots & a_1
\end{pmatrix}.\notag
\end{align}
\end{lemma}
\begin{proof}
From the hypothesis,
\begin{align}
\sum_{k=0}^{j} a_k c_{j-k} = b_j
\tag{1}
\end{align}
for each $j$.
Since $c_0 = b_0/a_0$,
\[
a_0 c_j + a_1 c_{j-1} + \cdots + a_{j-1} c_1=
(a_0 b_j - a_j b_0)/a_0,
\notag
\]
so, for $n \geq 1$,
\begin{align}
\begin{pmatrix}
a_0 & 0 & 0 & \cdots & 0 \\
a_1 & a_0 & 0 & \cdots & 0 \\
a_2 & a_1 & a_0 & \cdots & 0 \\
\vdots & \vdots & \vdots & \ddots & \vdots \\
a_{n-1} & a_{n-2} & a_{n-3} & \cdots & a_0
\end{pmatrix}
\begin{pmatrix}
c_1 \\ c_2 \\ c_3 \\ \vdots \\ c_n
\end{pmatrix}
= 
-\begin{pmatrix}
 (a_1 b_0 - a_0 b_1)/a_0 \\
(a_2 b_0 - a_0 b_2)/a_0  \\
(a_3 b_0 - a_0 b_3)/a_0  \\
\vdots \\
(a_n b_0 - a_0 b_n)/a_0 
\end{pmatrix}.
\tag{2}
\end{align}
Let 
\begin{align}
A_n =
\begin{pmatrix}
a_0 & 0 & 0 & \cdots & 0 \\
a_1 & a_0 & 0 & \cdots & 0 \\
a_2 & a_1 & a_0 & \cdots & 0 \\
\vdots & \vdots & \vdots & \ddots & \vdots \\
a_{n-1} & a_{n-2} & a_{n-3} & \cdots & a_0
\end{pmatrix}
\tag{3}
\end{align}
and let
\begin{align}
B_n = 
\begin{pmatrix}
a_0 & 0 & 0 & \cdots & -(a_1 b_0-a_0 b_1 )/a_0  \\
a_1 & a_0 & 0 & \cdots & -(a_2 b_0-a_0 b_2)/a_0  \\
a_2 & a_1 & a_0 & \cdots & -(a_3 b_0-a_0 b_3)/a_0  \\
\vdots & \vdots & \vdots & \ddots & \vdots \\
a_{n-1} & a_{n-2} & a_{n-3} & \cdots & -(a_n b_0-a_0 b_n)/a_0 
\end{pmatrix},\tag{4}
\end{align}
so that
\begin{align}
c_n = \frac{|B_n|}{|A_n|} =  \frac{|B_n|}{a_0^n}.
\tag{5}
\end{align}
Let
\begin{align}
D_n = \begin{pmatrix}
-(a_1 b_0 - a_0 b_1)/a_0 & a_0 & 0 & \cdots & 0 \\
-(a_2 b_0 - a_0 b_2)/a_0 & a_1 & a_0 & \cdots & 0 \\
-(a_3 b_0 - a_0 b_3)/a_0 & a_2 & a_1 & \cdots & 0 \\
\vdots & \vdots & \vdots & \ddots & \vdots \\
-(a_n b_0 - a_0 b_n)/a_0 & a_{n-1} & a_{n-2} & \cdots & a_1
\end{pmatrix},\notag
\end{align}
so $|D_n|=-|M_n|/a_0$ and $|B_n|= (-1)^{n-1} |D_n|$. Finally,
\begin{align}
c_n=
\frac{|B_n|}{a_0^n} =  \frac{(-1)^{n-1} |D_n|}{a_0^n}
=  \frac{(-1)^{n-1} (-|M_n|/a_0)}{a_0^n}
=\frac{(-1)^n}{a_0^{n+1}} \left |M_n \right |.
\tag{6}
\end{align}
\end{proof}
\begin{lemma}
The equations  below are equivalent. (We will
refer to both of them as equation (D) in the sequel.)
Let $h_0=1$ and
\[
n h_n = \sum_{r=1}^{n} j_r h_{n-r}\tag{D1}
\]
for $n \ge 1$.
With
$
H(t) = \sum_{n=0}^\infty h_n t^n$ and  $J(t) = \sum_{r=1}^\infty j_r t^r$,
\[
t \frac{d}{dt}H(t) = H(t) J(t).\tag{D2}\]
\end{lemma}
\underline{Claim}: Equation (D1) implies 
equation (D2).
\begin{proof}
From equation (D1),
$$
t \frac{d}{dt}H(t)  =
\sum_{n \ge 1}n h_nt^n = \sum_{n \ge 1}\left(\sum_{r=1}^{n} j_r h_{n-r}\right)t^n,
$$

$$
\text{while }
H(t) J(t) = 
\left( \sum_{k=0}^\infty h_k t^k\right)
\left(\sum_{r=1}^\infty j_r t^r \right)
=
\sum_{n=1}^{\infty}
\left
(
\sum_{\substack{k=n-r\\
k \ge 0 \\r\ge 1}} h_k  j_r
\right
)t^n =
$$

$$
\sum_{n=1}^{\infty}
\left
(
\sum_{\substack{
n-r \ge 0 \\ \\r\ge 1}} h_{n-r}  j_r
\right
)t^n =
\sum_{n=1}^{\infty}
\left
(
\sum_{
1 \le r \le n } h_{n-r}  j_r
\right
)t^n.
$$
\end{proof}
\underline{Claim}: 
(D2) implies (D1).
\begin{proof}
Re-ordering some of the equations
from the proof of the first claim,
we have
$$
\sum_{n \ge 1}n h_nt^n = 
t \frac{d}{dt}H(t)  = H(t) J(t)
=\sum_{n=1}^{\infty}
\left
(
\sum_{1 \le r \le n} h_{n-r}  j_r
\right
)t^n,
$$
and the claim follows by
equating coefficients.
\end{proof}
We introduce notation from the theory of partitions, mostly following the set-up in MacDonald's book ~\cite{M}, where a partition $\lambda$ is a non-increasing, possibly infinite sequence $(\lambda_1, \lambda_2,\dots)$  of non-negative integers containing only 
finitely many non-zero terms.
We write 
$|\lambda| = \sum_i \lambda_i$; $l(\lambda)$ for the number of positive parts in $\lambda$; $p(n) = p_n =
\sum_{|\lambda|=n} 1$;
and we set $z_{\lambda} = \prod_k k^{m_k(\lambda)} m_k(\lambda)!$. For sequences $\{s_0, s_1, ... \}$ such that $s_0 = 1$, we write
     \(s_\lambda =\prod_i s_{\pi_i}\). 
With \(f_0 = 1\),    
let \(\overline{f}\) denote the sequence 
\(f_0, f_1, f_2, ...\).
Let $J_n(\overline{j})$ and $H_n(\overline{h})$ be the matrices 
\[ 
\left(
\begin{matrix}
j_1 & -1 & 0 & \cdots & 0 \\
j_2 & j_1 & -2 & \cdots & 0 \\
\vdots & \vdots & \vdots & \vdots & \vdots \\
j_{n-1} & j_{n-2} &  \cdots & j_1 & -n + 1\\
j_n & j_{n-1} & j_{n-2} & \cdots & j_1\\
\end{matrix}
\right )
\]
and
\[ 
\left(
\begin{matrix}
h_1 & 1 & 0 & \cdots & 0 \\
2h_2 & h_1 & 1 & \cdots & 0 \\
\vdots & \vdots & \vdots & \vdots & \vdots \\
nh_n & h_{n-1} & h_{n-2} & \cdots & h_1\\
\end{matrix}
\right),
\]
respectively.
Equations (7) and (8)  in the following lemma appear in  Turnbull \cite{T}, Gould 
 ~\cite{Go99} 
and MacDonald ~\cite{M},
but the proofs are absent. 
We have not found suitable references to the original arguments.
\begin{lemma}
Let $h_0= j_1 = 1$.
    Equation (D) implies
that
\[
j_n = (-1)^{n+1} | H_n(\overline{h})|. \tag{7}
\]
and
 \[
n! h_n =  | J_n(\overline{j})|\tag{8}
\]
for $n \ge 1$.
\end{lemma}
\begin{proof}
First we derive equation (7).
From equation (D),
$J(t) = t H'(t)/H(t)$, so 
Lemma 2.1 tells us
that
\begin{align}
j_n = \frac{(-1)^n}{a_0^{n+1}} \left |M_n \right |,\notag
\end{align}
where $M_n$ is the $n \times n$ matrix
\[
\begin{pmatrix}
a_1 b_0 - a_0 b_1 & a_0 & 0 & \cdots & 0 \\
a_2 b_0 - a_0 b_2 & a_1 & a_0 & \cdots & 0 \\
a_3 b_0 - a_0 b_3 & a_2 & a_1 & \cdots & 0 \\
\vdots & \vdots & \vdots & \ddots & \vdots \\
a_n b_0 - a_0 b_n & a_{n-1} & a_{n-2} & \cdots & a_1
\end{pmatrix},
\]
$$
\sum_{k=0}^{\infty} b_k t^k = t H'(t) = 
  \sum_{k=0}^{\infty} k h_k t^k,
$$ and
$$
\sum_{k=0}^{\infty} a_k t^k
= 
H(t) =
\sum_{k=0}^{\infty} h_k t^k.
$$
Now $a_0 = 1$ and $b_0 = 0$, so 
\[M_n =
\begin{pmatrix}
 -  h_1 & h_0 & 0 & \cdots & 0 \\
  -  2h_2 & h_1 & h_1 & \cdots & 0 \\
  -  3h_3 & h_2 & h_2 & \cdots & 0 \\
\vdots & \vdots & \vdots & \ddots & \vdots \\
 -  n h_n & h_{n-1} & h_{n-2} & \cdots & h_1
\end{pmatrix},
\]
which yields equation (7).

To derive equation (8),
we rewrite equation (D1) 
as a recursion for $j_k$
\[
\sum_{r=1}^{k-1} -j_r h_{k-r} + k h_k = j_k,
\tag{D2}
\]
and summarize equation (D2) for $1 \le k \le n$ as 
\[
\begin{pmatrix}
1 & 0 & 0 & \dots & 0 \\
-j_1 & 2 & 0 & \dots & 0 \\
-j_2 & -j_1 & 3 & \dots & 0 \\
\vdots & \vdots & \vdots & \ddots & \vdots \\
-j_{n-1} & -j_{n-2} & -j_{n-3} & \dots & n
\end{pmatrix}
\begin{pmatrix} h_1 \\ h_2 \\ h_3 \\ \vdots \\ h_n \end{pmatrix} =
\begin{pmatrix} j_1 \\ j_2 \\ j_3 \\ \vdots \\ j_n \end{pmatrix}.
\tag {D3}
\]
Let $A$ be the left-most matrix in equation (D3),
so that $|A| = n!$, and let
$$
B = \begin{pmatrix}
1 & 0 & 0 & \dots & 0 & j_1 \\
-j_1 & 2 & 0 & \dots & 0 & j_2 \\
-j_2 & -j_1 & 3 & \dots & 0 & j_3 \\
\vdots & \vdots & \vdots & \ddots & \vdots & \vdots \\
-j_{n-1} & -j_{n-2} & -j_{n-3} & \dots & -j_1 & j_n
\end{pmatrix}.
$$
Equation (D3) gives
$h_n = |B|/|A| = |B|/n!$.
But $|B| = |J_n(\overline{j})|$.
\end{proof}
The proof of Lemma 2.4 is essentially identical
to calculations from section I.2 of 
MacDonald's book ~\cite{M}.
\begin{lemma}\footnote{This may be considered a generalization of some examples later in this article, for example, Theorem 3.1(i). A. Petojevi\'{c} has shown the writer another 
generalization of that result. It is part of the draft article ``An Explicit Formula
of the Coefficients of Eta-Products'' (by 
Petojevi\'{c} and S. Orli\'{c}),
stored as ``petojevic's draft3oct25'', in the repository \cite{Bre25}.}
    Equation (D) implies
\[
h_n = \sum_{|\lambda| = n} z^{-1}_{\lambda} j_{\lambda}.
\tag{9}
\]
\end{lemma}
\begin{proof} 
Let
\(H_1(t) = H(t) -1\) and 
\(J_1(t) = J(t)/t\).
\rm
We have \[
t H_1(t) J_1(t) = 
\left( \sum_{s \geq 1} h_s t^s \right)
\left(\sum_{r \geq 1} j_r t^r\right)
= 
\sum_{n=1}^{\infty} \left( \sum_{r=1}^n  j_r h_{n-r}  \right) t^n
\]
and
\[ 
t H_1'(t)
= \sum_{n=1}^{\infty} n h_n t^n
= \sum_{n=1}^{\infty} \left(\sum_{r=1}^n  j_r h_{n - r} \right)t^n,
\notag
\]
so
\[
t H_1(t) J_1(t) = t H_1'(t).
\]
Therefore
\[
J_1(t) = H_1'(t)/H(t) = (\log H_1)',
\]
so
\[
H_1(t) = \exp \sum_{r \geq 1} (j_r t^r/r) = \prod_{r\geq 1} \exp(j_r t^r /r)
= \]\[\prod_{r \geq 1}\left(\sum_{m_r=0}^{\infty} 
\frac{(j_r t^r)^{m_r}}{r^{m_r}}
m_r !\right) =
\sum_{\lambda}  z_{\lambda}^{-1} j_{\lambda} 
t^{|\lambda|}
= \sum_{n\ge 1}\left(\sum_{|\lambda| = n} z^{-1}_{\lambda} j_{\lambda} t^n\right).
\]
\end{proof}
 \section{Examples}
 In the sequel, we
 will move freely back and forth between the notations 
 \(f(n)\) and \(f_n\). 
 Ramanujan's recursion for $\tau(\nu)$  (\cite{R}, equation (99), page 152)
\[
\tau(\nu)=\frac{24}{1-\nu}\sum_{k=1}^{\nu-1} \sigma(k)\tau(\nu-k) \tag{10}
\]
requires \(\nu > 1\).
Setting \(n = \nu-1\), we get
\[
n \tau(n+1) =\sum_{k=1}^{n} -24\sigma(k) \tau(n+1 - k). 
\tag{11}
\]
For comparison with later statements, we write $T_n = \tau(n+1), S_k = -24\sigma(k)$,
and get
\begin{equation}
n T_n = \sum_{k=1}^{n} S_k T_{n-k}.
\tag{11\textquotesingle{}}
\end{equation}
With $h_n = T_n = \tau(n+1)$ and $j_n = S_n=-24\sigma(n)$, the lemmas imply the next theorem.
\begin{theorem}
\[
n! T_n = | J_n(\overline{S})| \tag{i}.
\]
\[
S_n = (-1)^{n-1} | H_n(\overline{T})| \tag{ii}.
\]
\[
-24\sigma(n) = (-1)^{n-1} |H_n(\overline{T})|
\tag{iii}.
\]
\[
T_n = \sum_{|\lambda| = n} z^{-1}_{\lambda} S_{\lambda}
\tag{iv}.
\]
\end{theorem}
\begin{remark} 
We will use \(
\chi_M
\) to denote the characteristic polynomial of a matrix 
\(
M
\). In particular, we make the following choice of sign (which is not uniformly applied uniformly in the literature):
\(\chi_M = |xI_n - M| \).
Then the following statements are equivalent:  
\(
\tau(n+1) = 0
\), \(|J_n(\overline{S})| = 0\) and 
\(
\chi_{J_n(\overline{S})}(0)  = 0
\).
\end{remark}
For brevity, let us write \(
\Sigma_k = J_k(\overline{S}),
r_a = 1/a!\)
and \(m*n = (m+1)(n+1)\).
\begin{corollary}
If \(m+1, n+1\) are coprime positive integers, then
\[r_m r_n
|\Sigma_m|\cdot
|\Sigma_m|
=
r_{m*  n}|\Sigma_{m*n}|.\tag{12}
\]
\end{corollary}
\begin{proof}
\[  r_{m*n} |J_{m*n}(\overline{S})| =
T_{m*n}=
T_{m+1} T_{n+1} =
 r_{m+1} r_{n+1} | J_{m+1}(\overline{S})|
 | J_{n+1}(\overline{S})|.
\] 
\end{proof}
The corollary below began life as the present writer's conjecture, based on computer data. His colleague L\'opez-Bonilla supplied the proof below, which 
uses Theorem 4.23 in Vein and Dale's 
book ~\cite{Vein1999}. Aleksandar Petojevi\'{c} sent him another proof, which is on his repository for this article
~\cite{Pe}.
\begin{corollary}(L\'opez-Bonilla, \text{personal communication})
    \[
    \chi_{\Sigma_n}
    =
    \sum_{k=0}^{n} \frac{n!}{k!}\tau(n-k+1) x^k.
    \]
\end{corollary}
To state L\'opez-Bonilla's proof, we first reproduce Vein and Dale's theorem. 
Let 
\(A_n =\)
\[
\begin{pmatrix}
a_1 & -1 & 0 & \cdots & 0 \\
a_2 & a_1 & -2 & \cdots & 0 \\
\vdots & \vdots & \vdots & \ddots & \vdots \\
a_{n-1} & a_{n-2} & \cdots & a_1 & -(n{-}1) \\
a_n & a_{n-1} & \cdots &a_2& a_1
\end{pmatrix}^T,
\]
where the $T$ operator is matrix transpose.
Let \(B_n(x) = |A_n -xI_n|\),
so that 
\(B_n(x) = (-1)^n \chi_{A_n}(x)\). 
 Then 
\begin{proposition}
(Vein and Dale, Theorem 4.23)
\[
B_n(x) = \sum_{r=0}^n 
\binom{n}{r}
|A_r|(-x)^{n-r}.
\]
\end{proposition}
Here is Lopez-Bonilla's argument for Corollary 3.3:
\begin{proof}
Let \( r \) range from \( 1\)  to \( n \).
With \(a_k = -24\sigma(k), 
A_r^T = -\Sigma_r\).
Theorem 3.2 implies that
\[
|A_r|=|A_r^T|
= | -\Sigma_r| =
 (-1)^r|\Sigma_r| = r! \tau(r+1).\tag{13}
\]  
Let \(\lambda = -x\). Then
\[
\chi_{\Sigma_n}(x)= |xI_n -J_n(\overline{S})|
= |xI_n - (-A^T_n)|=\]
\[
\chi_{-A^T_n}(x)=\chi_{-A_n}(x)=
|x I_n - (-A_n)|=
\]\[
|-\lambda  I_n - (-A_n)|=
B_n(\lambda) =\]
(by Proposition 2)
\[
\sum_{r=0}^n 
\binom{n}{r}
|A_r| (-\lambda)^{n-r} = 
\]
(by equation (13))
\[
\sum_{r=0}^n 
\binom{n}{r}
 r! \tau(r+1) x^{n-r}=
\]
\[
\sum_{r=0}^n 
\frac{n!}{(n-r)!}
 \tau(r+1)x^{n-r}
= \]\[
\sum_{k=0}^n 
\frac{n!}{k!}
 \tau(n-k+1)x^k.
\]
\end{proof}
  Let $\chi_{H_n(\overline{T})}(x) = \sum_{k=0}^{n}  \gamma_{n,k} x^k$. 
  The sequence \(
\{c_n
\}\) \rm mentioned below
 has a page on the On-line Encyclopedia of Integer Sequences site ~ \cite{SloaneA006922}.
Our calculations are accessible on this article's
GitHub repository ~\cite{Bre25}.
On the basis of those calculations, we make the following conjecture.
\begin{conjecture}
1. $\gamma_{n,0} = 24 \sigma(n).$
2. Let \(\eta(z)\) be Dedekind's eta function and
let \(c_n\) be the coefficient of \(q^n\)
in the Fourier expansion of $1/\eta(q)^{24}$
(its expansion in powers of \(\exp 2 \pi i z\)).
Then \(
\gamma_{n,1}= (n+1) c_n
\).
% see cn1_11aug25.ipynb
\end{conjecture}
\begin{remark}
We notice, of course, that the matrix $H_n(\overline{T})$
is itself populated with Fourier coefficients of $\Delta = q \cdot\eta(q)^{24}$.
\end{remark}
We recall the  identity
~\cite{HR1918} (or page 279 in Ramanujan's collected papers ~\cite{R2015})
\[
n p(n) = \sum_{k=1}^n \sigma(k)p(n-k).\tag{14}
\]
Evidently, Hardy and Ramanujan take 
\(p(0)\) equal to zero in this identity.
Therefore  equation (1) is satisfied with 
\(\overline{h} = \overline{p}\) 
and \(\overline{j} = \overline{\sigma}\).
So the lemmas produce
\begin{theorem}
    \[
n! p_n =  | J_n(\overline{\sigma})|,\tag{i}
\]
\[
\sigma_n = (-1)^{n-1} | H_n(\overline{p})|,\tag{ii}
\]
and
\[
p_n = \sum_{|\lambda| = n} z^{-1}_{\lambda} \sigma_{\lambda}.
\tag{iii}
\]
\end{theorem}
\begin{remark}
The notation  ``\(\sigma_{\lambda}\)'' an example of the usage borrowed 
from MacDonald) and explained in section 2 above, with
\(s_n = \sigma_n.\)
\end{remark}
By Theorem 3.4, we can substitute 
\(\frac{(-1)^k}{24}|H_k(\overline{T})|\) 
for \(\sigma(k)\) in equation (16),
and then we have
\begin{theorem}
\[
n p_n = \sum_{k=1}^n \frac{(-1)^k}{24}|H_k(\overline{T})| p(n-k).\tag{15}
\]
\end{theorem}
This is a version of equation (1) in which
\(h_n = p_n\),  \(h_{n-r} = p_{n-r}\), and \(j_r =  \frac{(-1)^r}{24}|H_r(\overline{T})| = \phi_r\) (say). 
Now the lemmas produce the following. 
\begin{theorem}
    \[
n! p_n =  | J_n(\overline{\phi})|,\tag{i}
\]
\[
\phi_n = (-1)^{n-1} | H_n(\overline{p})|,\tag{ii}
\]
and
\[
p_n = \sum_{|\lambda| = n} z^{-1}_{\lambda} \phi_{\lambda}.
\tag{iii}
\]
\end{theorem}
\section{Deformations}
\subsection{Definition}
We define 
deformations $J^{(c)}_n(\overline{j})$ of the $J_n(\overline{j})$ as matrices
\[ 
\left(
\begin{matrix}
c & -1 & 0 & \cdots & 0 \\
j_1 & c & -2 & \cdots & 0 \\
\vdots & \vdots & \vdots & \vdots & \vdots \\
j_{n-2} & j_{n-3} &  \cdots & c & -n + 1\\
j_{n-1} & j_{n-2} & j_{n-3} & \cdots & c\\
\end{matrix}
\right).
\]
Our original reason for focusing on them was a phenomenon we will discuss in the section about characterisric polynomials following this one;
we came across it because
of a fortuitous coding mistake.
Here we describe a combinatorics-flavored
observation about these objects.

\subsection{Derangements}\footnote{See ``deformed one fn7oct25.ipynb'' in the \cite{Bre25}.}
Let $\iota(n) = 1 = h_{\iota}(n)$ (say) for all non-negative integers $n$. Let $\overline{j_{\iota}}$ be the companion sequence to
$\overline{h_{\iota}}$ in the sense of equation (D).
In the parlance of combinatorics, a permutation of $n$ 
objects with no fixed point is a \emph{derangement}, and the number of such permutations is represented as $!n$. 
\begin{conjecture}
$|J_n^{(1)}(\overline{j})|= !n$.
\end{conjecture}
\noindent
We verified this 
for $1 \le n \le 300$. 

There are two  recursions by Euler for $!n$ on page A000166 in the Encyclopedia
of Integer Sequences ~\cite{OEISA000166}. Our guess is that recursion (D), in this case, is  equivalent to the following generating function
with $k=y=0$
$$
\sum_{k,n \ge 0}
D_k(n) x^n y^k /n! = 
\frac{e^{-x(1-y)}}{1-x},
$$
where $D_k(n)$ is the number of permutations
of $n$ letters with exactly $k$ fixed 
points (\cite{Wilf2006}, 2.7, exercise 27(d)), so that $!n = D_0(n)$.
\section{Geometry of the roots of the characteristic polynomials} 
\subsection{Characteristic polynomials of matrices representing the values tau(n th prime)
}
Because the smallest
$n$ such that $\tau(n) = 0$ 
(if such an $n$ exists)\footnote{
We discussed this question in section 9 of our 2021 article
\cite{Bre24}.} 
occurs at $n$ a prime number ~\cite{Le47},
we applied Lemma 2.3 to $h_n = 
\tau(p_n)$.
We
used equation (D1) to calculate initial segments of the companion sequence \(
\{j_n\}_{n \geq 1}\) corresponding to 
$\overline{h}=\{1,\tau(p_1), \tau(p_2),...\}$.
(Here $h_0 = 1$ as required in our lemmas.)
If none of the roots of the associated characteristic  polynomials $\Pi_n(x)$ (say) $=\chi_{J_n(\overline{j})}$ is ever zero, then neither the determinant $|\Pi_n|$  nor the corresponding value of $\tau(p_n)$ ever vanishes. For each $n$, we attempted to compute the minimum modulus of the roots of 
the associated characteristic  polynomials $\chi_{J_n(\overline{j})}$. 

Because of the aforementioned coding mistake. we first ran through this procedure with  deformations $J^{(1)}_n(\overline{j})$
and observed a phenomenon that 
is apparently
absent from the behavior of the corresponding 
undeformed objects. 

Let is write $\Pi^{(c)}_n = \chi_{J^{(c)}_n(\overline{j})}$. In Figure 1, we plot the minimum modulus among the roots of 
$\Pi^{(1)}_n(x)$ against $n$.
 Evidently, this graph is roughly periodic, and its upper- and lower-envelopes are roughly sinusoidal.\footnote{See the notebook ``deformations for sinusoidality 5 oct25no1.ipynb'', which is in our repository ~\cite{Bre25}, for Figures 1, 2 and 3.} 
\begin{figure}[htbp]
    \centering
    \hspace*{0em}
    \includegraphics[width=1\textwidth]{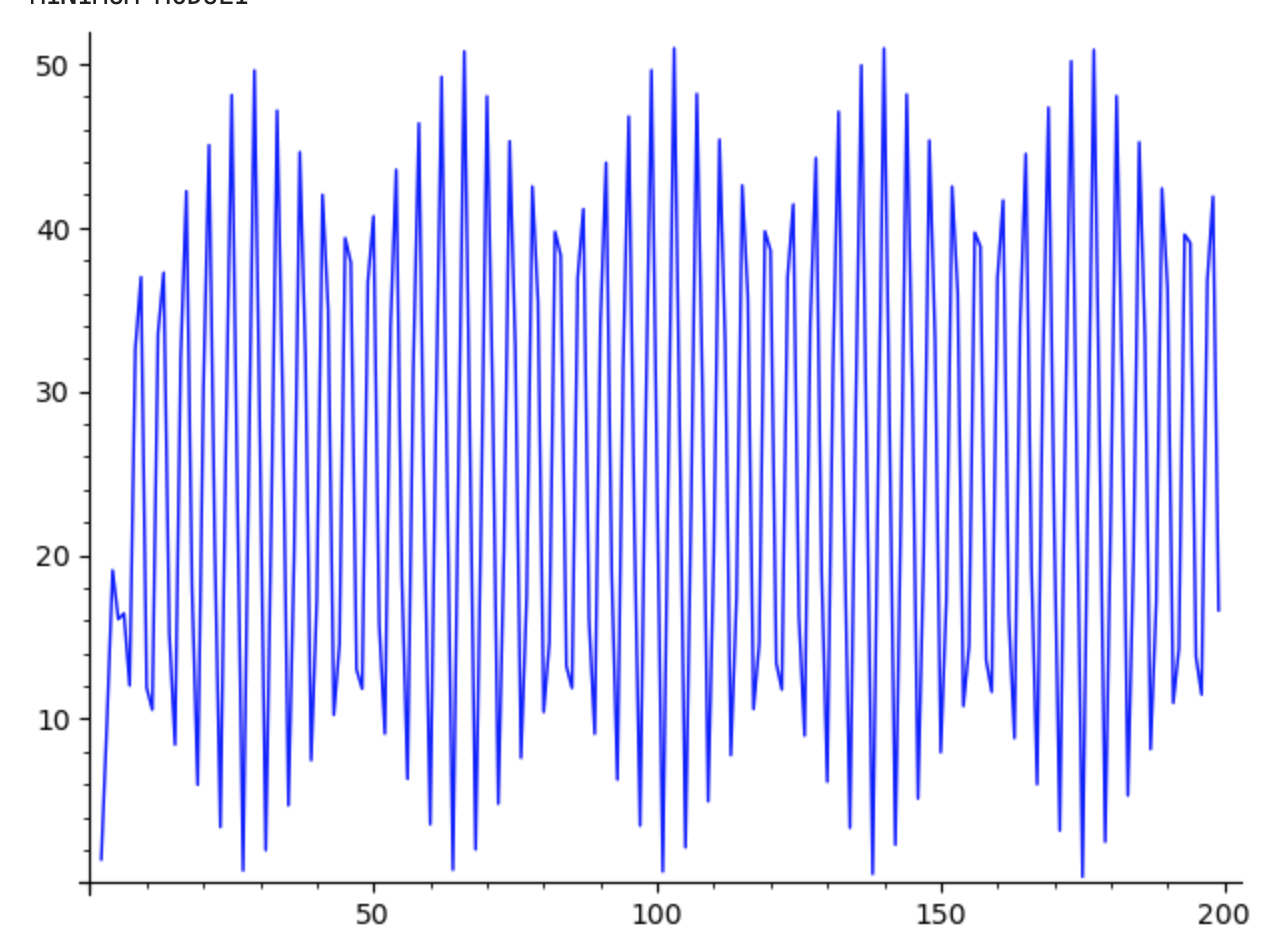}
    \caption{Minimum moduli of the roots of the
    \(
    \chi_{\Pi^{(1)}_n}
    \)} (deformed matrices with $c=1$.)
    \label{fig:new_primeTau_mins.png}
\end{figure}

We were initially under the impression that we were plotting the minima for the $\Pi_n$, so we wanted to make quite visible their failure to touch the $n$-axis. For that reason, in Figure 2, we plotted the logarithms of the minima. Because of our mistake, this plot is  irrelevant to Lehmer's question, but the lower envelope of the plot does
appear to suggest that the minima of the zeros of the $\chi_{\Pi_n^{(1)}}$ may have zero as a limit point.

To test the hypothesis that the entire zero-set of $\chi_{\Pi_n^{(1)}}$ is, in some loose sense, oscillating with $n$, we plotted the maxima 
as well (Figure 3.)  The graph is nearly a straight line with slope $24$, an interesting number under the circumstances. It seems to cast doubt on the hypothesis that
the root sets oscillate when regarded as geometric ensembles,
but we hope to study the question further.
\begin{figure}[htbp]
    \centering
    \hspace*{0em}
    \includegraphics[width=1\textwidth]{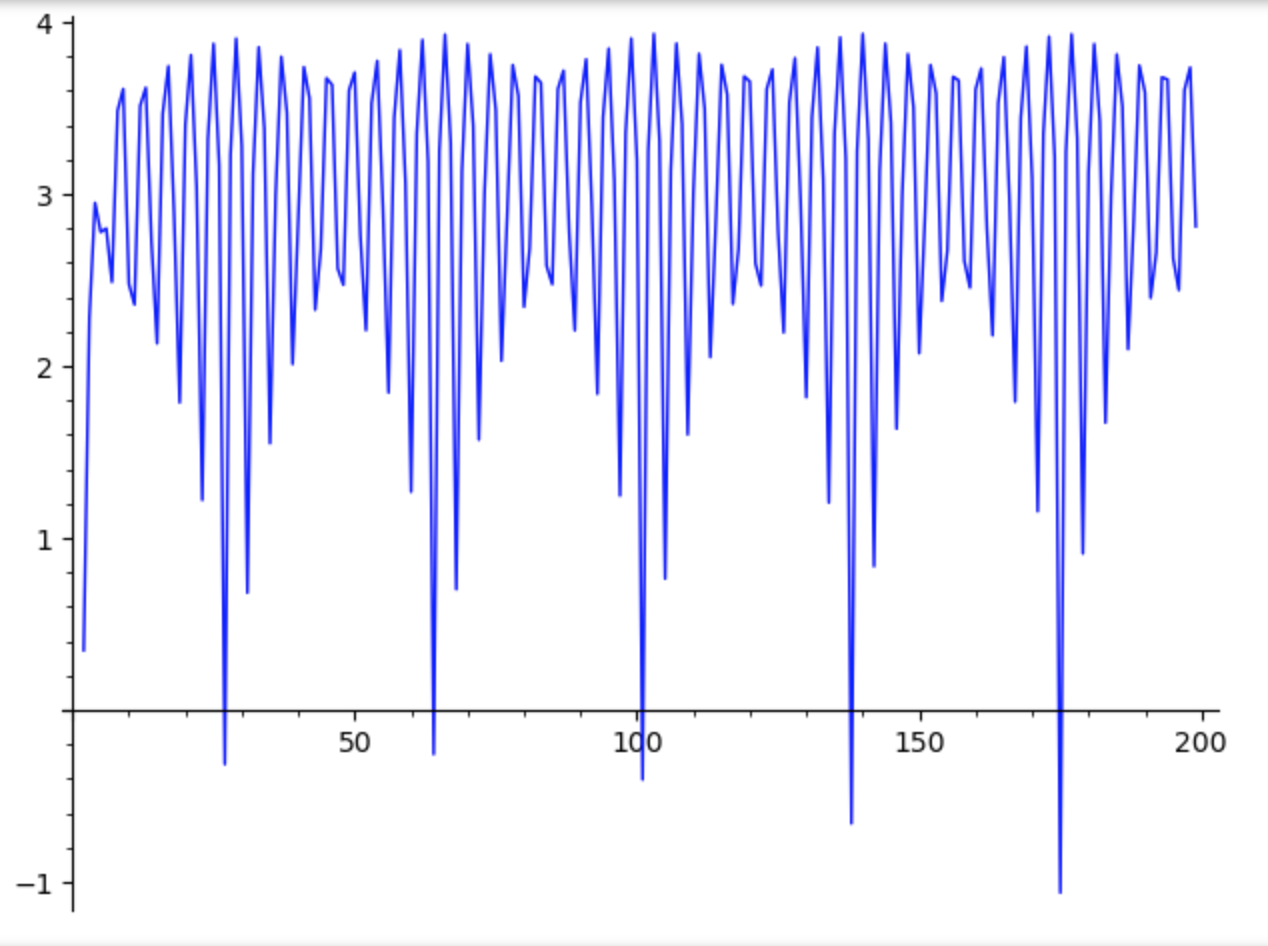}
    \caption{Logarithms of the minimum moduli of the
    roots of the
    \(
    \chi_{\Pi^{(1)}_n}
    \) (deformed matrices with $c=1$.)}
    \label{fig:new_primeTau_log_mins.png}
\end{figure}
\begin{figure}[htbp]
    \centering
    \hspace*{0em}
    \includegraphics[width=1\textwidth]{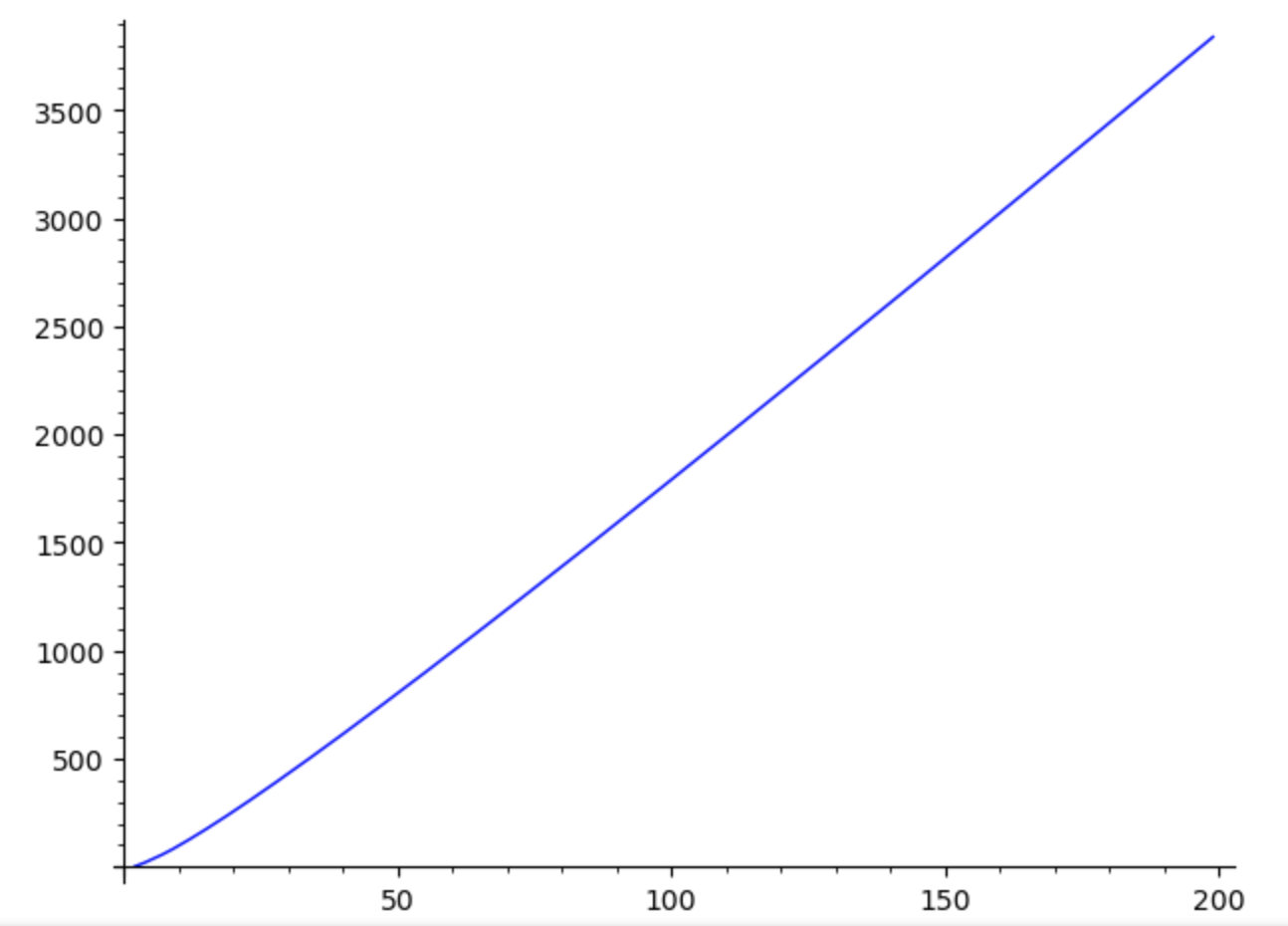}
    \caption{Maximum moduli of the roots of $\chi_{\Pi^{(1)}_n}$ (deformed matrices with $c=1$.)}
    \label{fig:new_primeTau_maxes.png}
\end{figure}

Other values of $c$ give similar behavior, but it was most pronounced among the values we tested at $c=1$. In particular, the roughly periodic oscillatory pattern is sometimes visible, but the 
sinusoidality of the upper- and lower- envelopes is not.

We of course must worry that the observed regularities have been produced in error. Or it may be that the regularities are a widespread effect, not associated with the underlying sequences, but only with the machinery of Lemma 2.3 more generally. We tested this by cooking up random sequences $\overline{j}$ and making the appropriate plots, but nothing like periodicity shows up.\footnote{See ``random sequence 28sept25.ipynb'' in ~\cite{Bre25}.}

Figure 4 shows the altogether irregular behavior
of the minimum moduli for the undeformed $\chi_{\Pi_n}$.\footnote{For Figures 8-10, see
``undeformed primeTau 5oct25.ipynb'' in ~\cite{Bre25}.} Figure 5 shows that the maximum moduli for these polynomials line up on a nearly straight line, and Figure 6 plots the slopes
of the chords from the points in Figure 5 to $(0,0)$. The curve appears to be asymptotic to a horizontal line with height approximately $22.924$.
\begin{figure}[htbp]
    \centering
    \hspace*{0em}
    \includegraphics[width=1\textwidth]{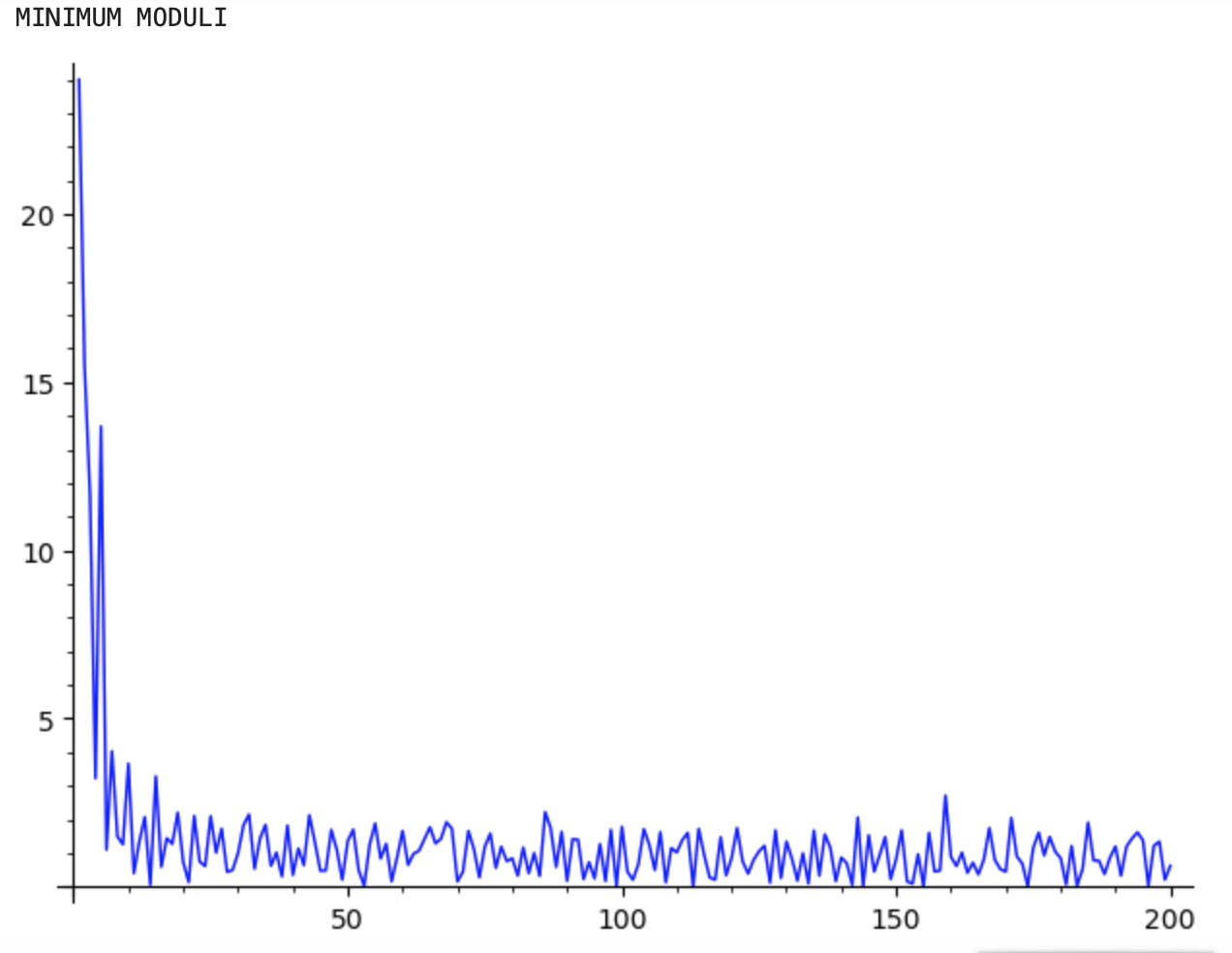}
    \caption{Minimum moduli of the roots of 
    undeformed
    $
    \chi_{\Pi_n}
    $}
    \label{fig:undeformed_primeTau_min_moduli.png}
\end{figure}
\begin{figure}[htbp]
    \centering
    \hspace*{0em}
    \includegraphics[width=1\textwidth]{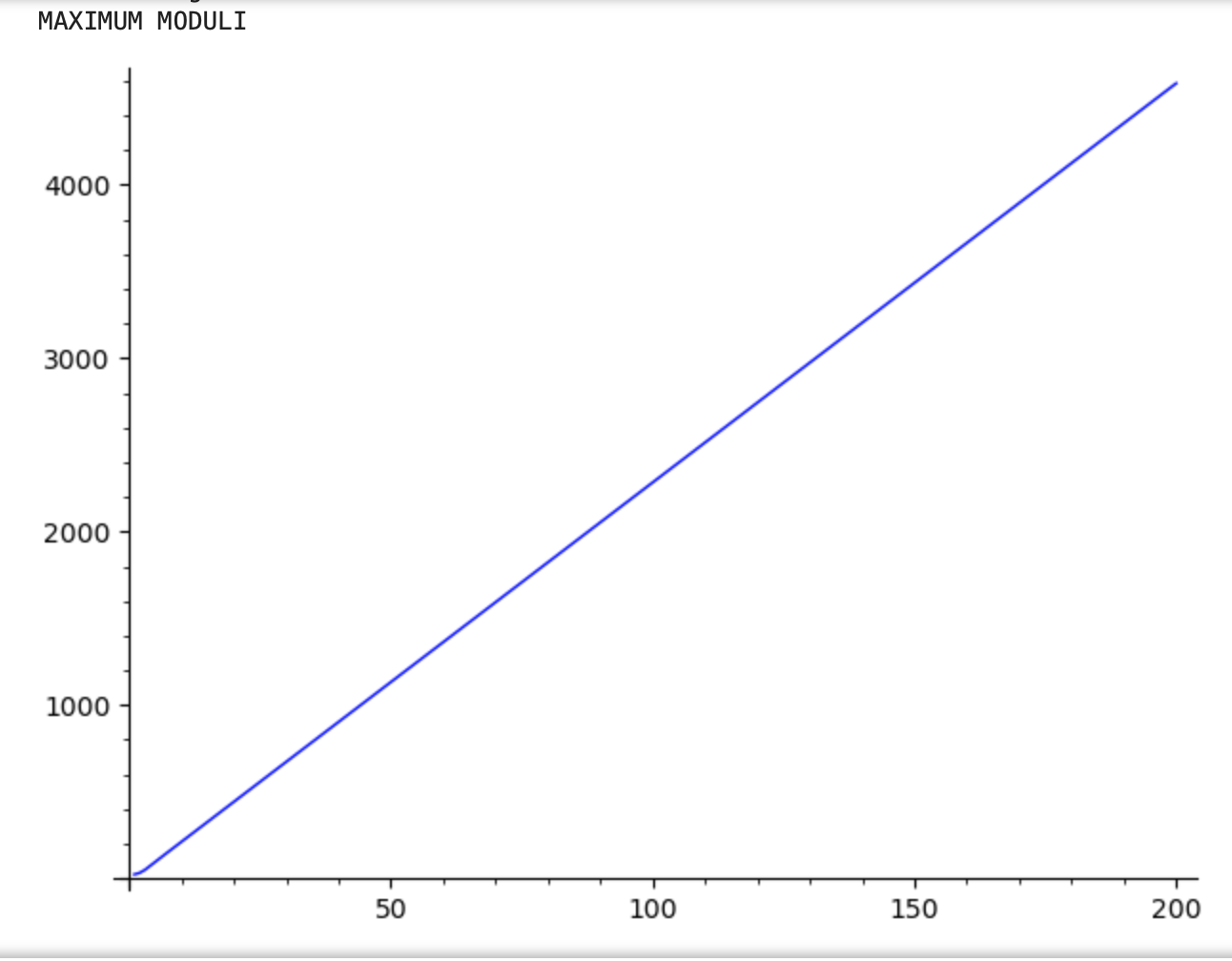}
    \caption{Maximum moduli of the roots of 
    undeformed
    $
    \chi_{\Pi_n}
    $}
    \label{fig:undeformed_primeTau_max_moduli.png}
\end{figure}
\begin{figure}[htbp]
    \centering
    \hspace*{0em}
    \includegraphics[width=1\textwidth]{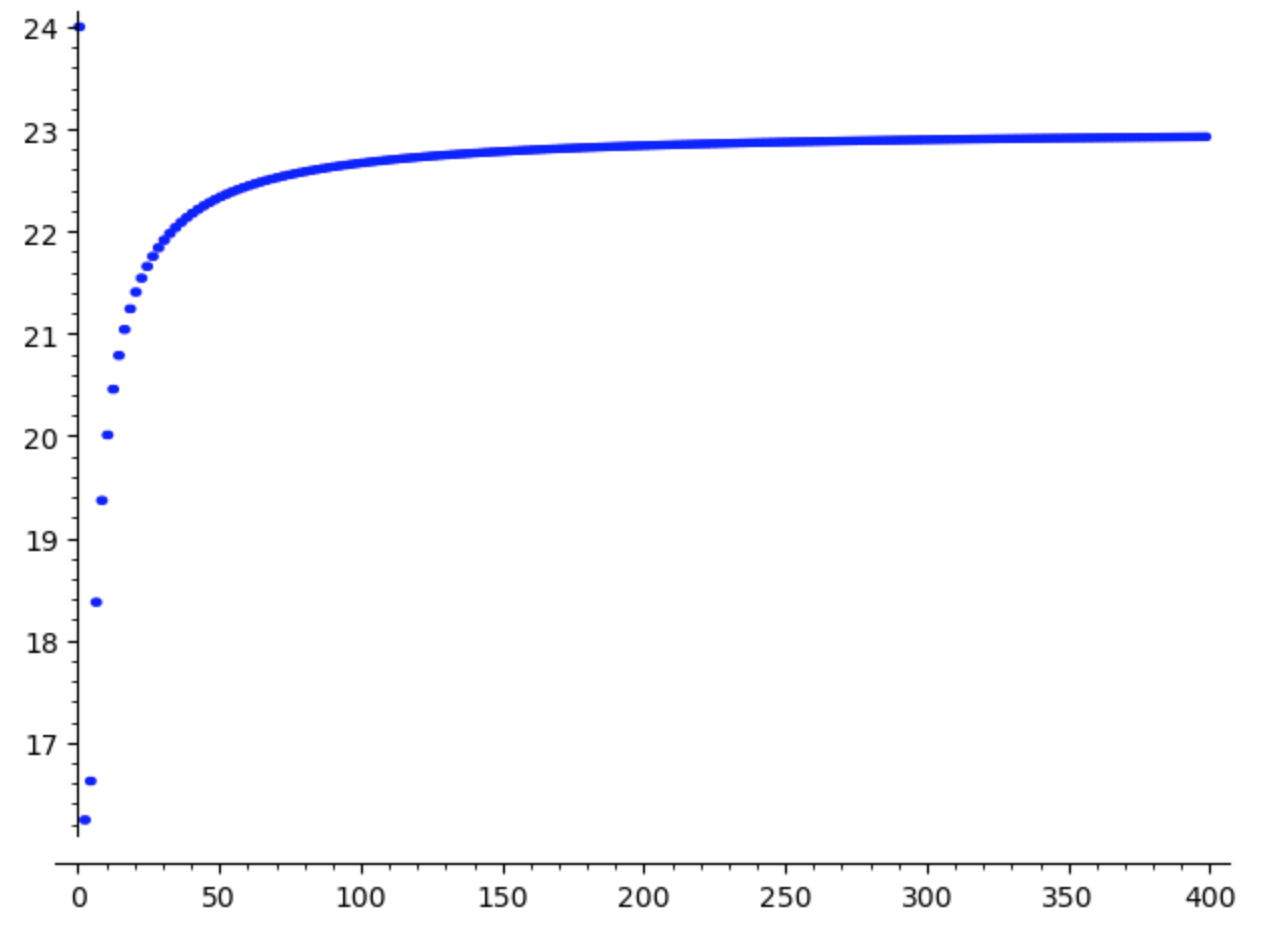}
    \caption{Slopes of partial plots of maximum moduli
    of roots of undeformed
    $
    \chi_{\Pi_n}
    $}
    \label{fig:slopes_undeformed_max_moduli.png}
\end{figure}
\subsection{Systematic searches for periodic behavior}\footnote{See 
``deformed pp minus c 8oct25.ipynb'' in
~\cite{Bre25}.}
Graphs such as that in Figure 1 led us to wonder
whether or not the lower envelopes of plots of the minimum moduli of roots of characteristic functions for matrices 
$J_n^{(c)}(\overline{j})$ -- where 
$j$ is the companion function 
in the sense of equation (D)
of a function
$h(n) = f(x)-v$ for certain $f$ and $v$ -- are sinusoidal.
So far, statistical tests fail to verify this,
 but we have had better luck testing simply for periodicity of the entire plots, using auto-corellation scores from bespoke code made by the AI Claude.

For example, let $P_v(n) = p_{n+1}-p_n - v$. The prime pair hypothesis is just the claim that the zero set of $P_2$ is infinite. Making plots seems unlikely to be a fruitful approach to this problem, of course, but we made ``deformed'' plots (in the sense above) of the relevant minimum moduli for several values of $v$ anyway and saw nearly-linear graphs, except for $v = -2, -1, 0, 1$ and $2$. We display those graphs in Figure 7.

Another example we found was the class of vertical translates $\mu_v(n) = \mu(n) - v$ for the M\"obius function.\footnote{See ``deformed mobius$_v$ 12oct25.ipynb'' in \cite{Bre25}. (We have altered our notation here to avoid a clash with our use of $c$ in the notation for deformed matrices introduced above, but the confusing use of $c$ is
retained, unfortunately, in our Jupyter notebooks.)} We searched near $v=0$ and found at increments of $0.2$, that
$v=0$ appears to be close to optimal.\footnote{We gathered the statistical reports in a folder of screenshots ``moebius$_c$ stats'' in \cite{Bre25}. (Same remarks as to the use of $v$ vs. $c$.)} As $v$ approaches zero, plotting the minimum modulus 
of the characteristic functions of the ``deformed'' matrices (with $c=1$) of size $n$ against $n$ itself,  there appears to be something like a phase transition, from an approximately straight line through the origin, to the approximately periodic behavior at $v=0$ (Figure 8.)
\begin{figure}[htbp]
    \centering
    \begin{subfigure}{0.45\textwidth}
        \centering
        \includegraphics[width=\textwidth]{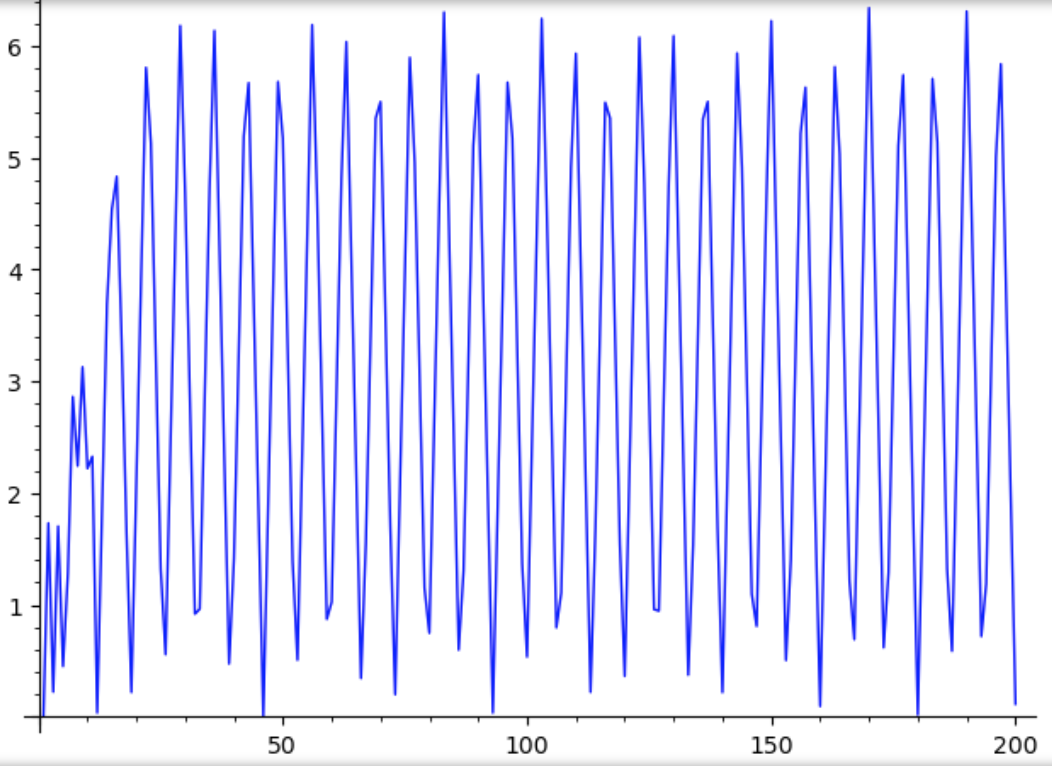}
        \caption*{$v =-2$}
        \label{fig:v=-2}
    \end{subfigure}
    \hfill
    \begin{subfigure}{0.45\textwidth}
        \centering
        \includegraphics[width=\textwidth]{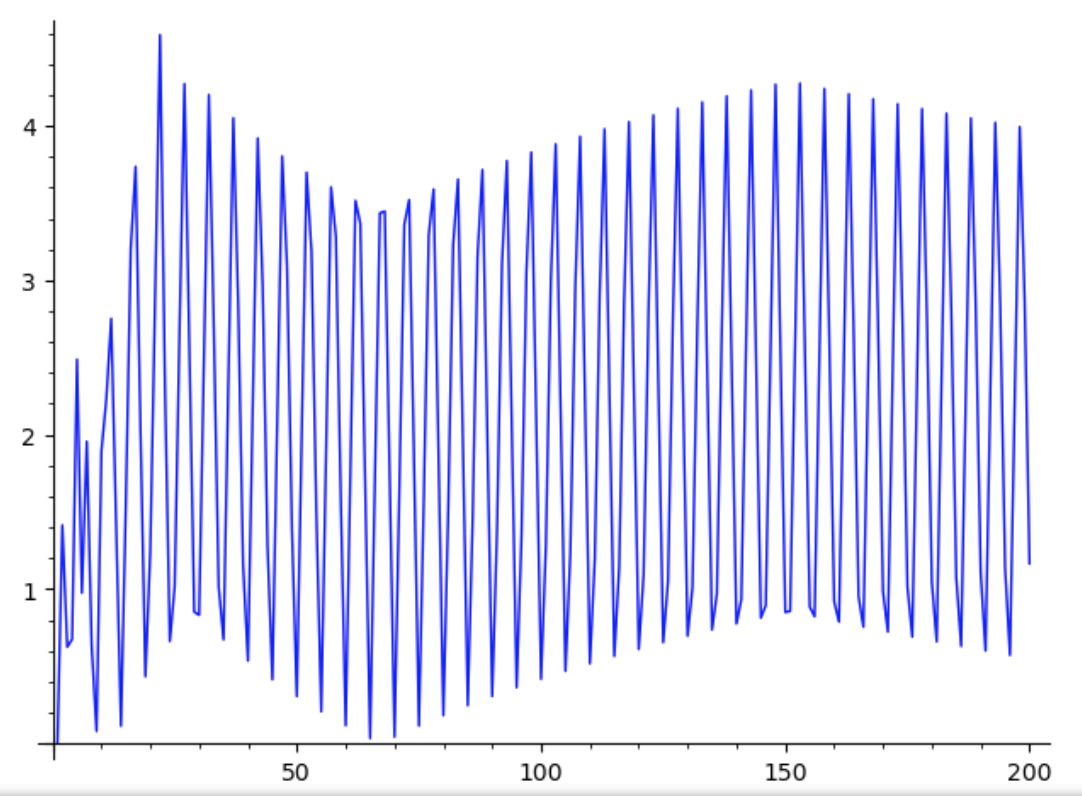}
        \caption*{$v = -1$}
        \label{fig:v=-1}
    \end{subfigure}
    
    \begin{subfigure}{0.45\textwidth}
        \centering
        \includegraphics[width=\textwidth]{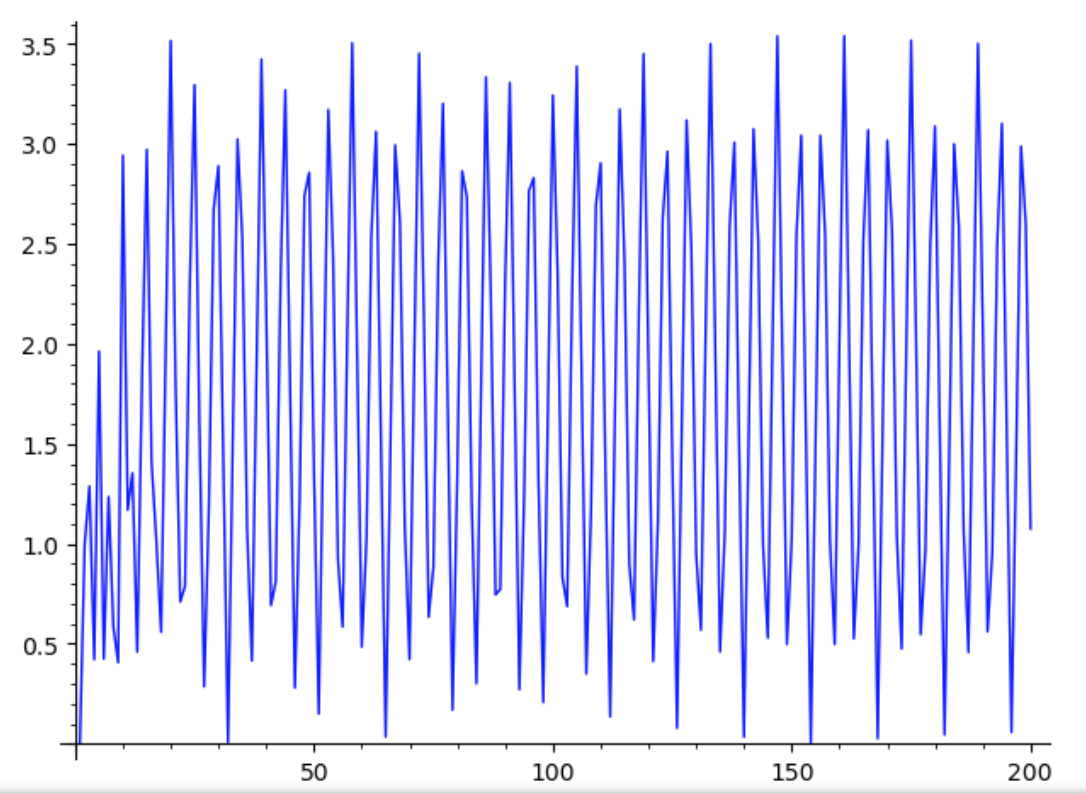}
        \caption*{$v = 0$}
        \label{fig:v=0}
    \end{subfigure}

    \vspace{1em} % Vertical space between rows
    \hfill
    \begin{subfigure}{0.45\textwidth}
        \centering
        \includegraphics[width=\textwidth]{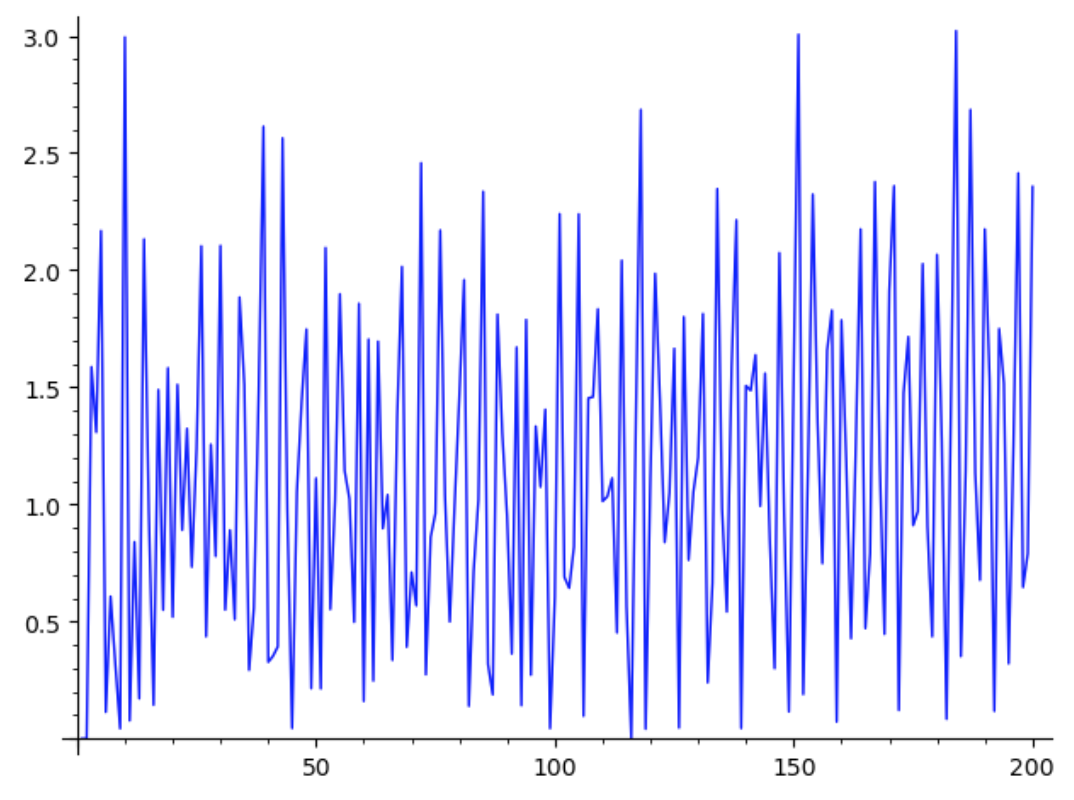}
        \caption*{$v=1$}
        \label{fig:v=1}
    \end{subfigure}
    \begin{subfigure}{0.45\textwidth}
        \centering
        \includegraphics[width=\textwidth]{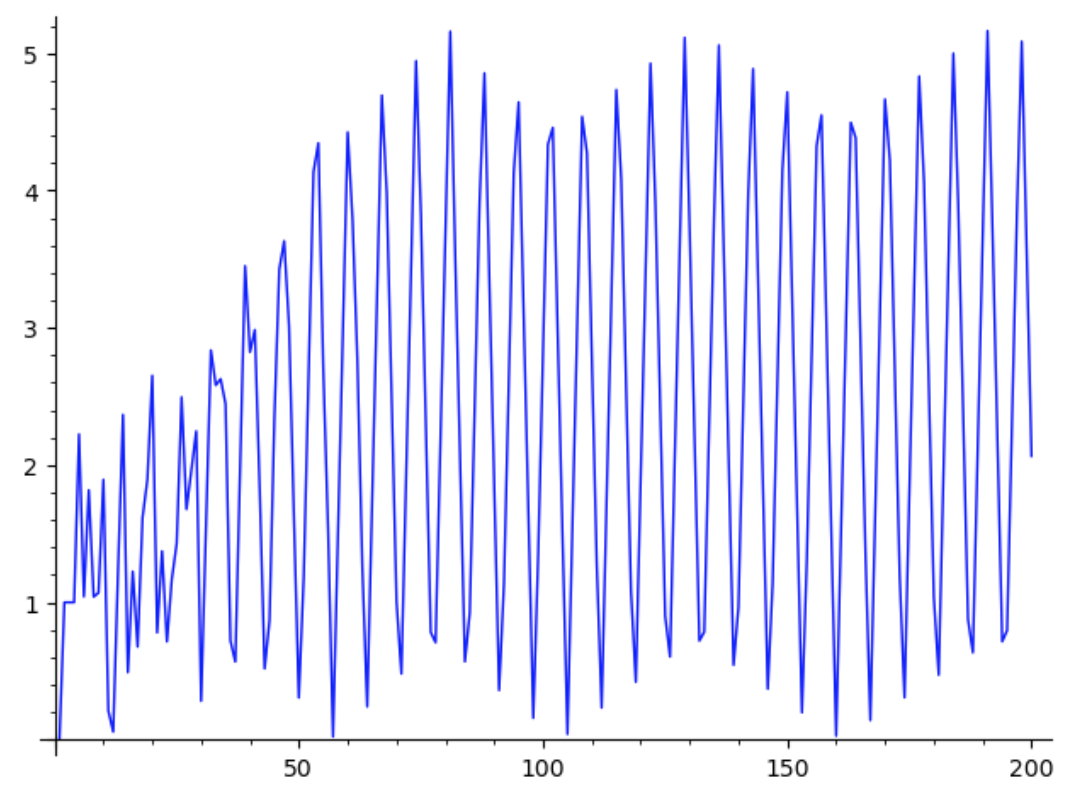}
        \caption*{$v=2$}
        \label{fig:v=2}
    \end{subfigure}
    
    \caption{Minimum moduli for roots of characteristic functions of matrices $J_n^{(1)}(\overline{j})$ associated by equation (D) to the function $h(n) = P_v(n)=p_{n+1}-p_n-v$.}
    \label{fig:main}
\end{figure}
\begin{figure}[htbp]
    \centering
    \begin{subfigure}{0.45\textwidth}
        \centering
        \includegraphics[width=\textwidth]{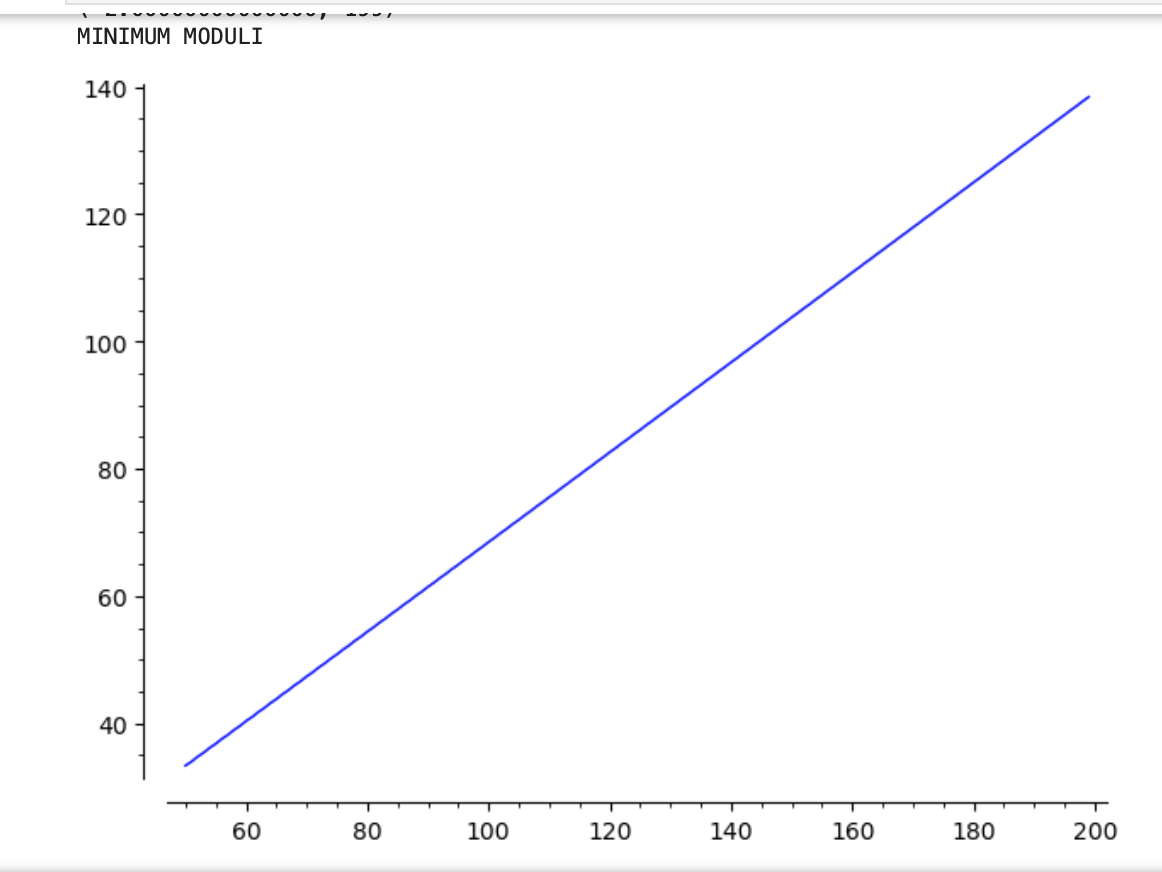}
        \caption*{$v =-2$}
        \label{fig:v=-2}
    \end{subfigure}
    \hfill
    \begin{subfigure}{0.45\textwidth}
        \centering
        \includegraphics[width=\textwidth]{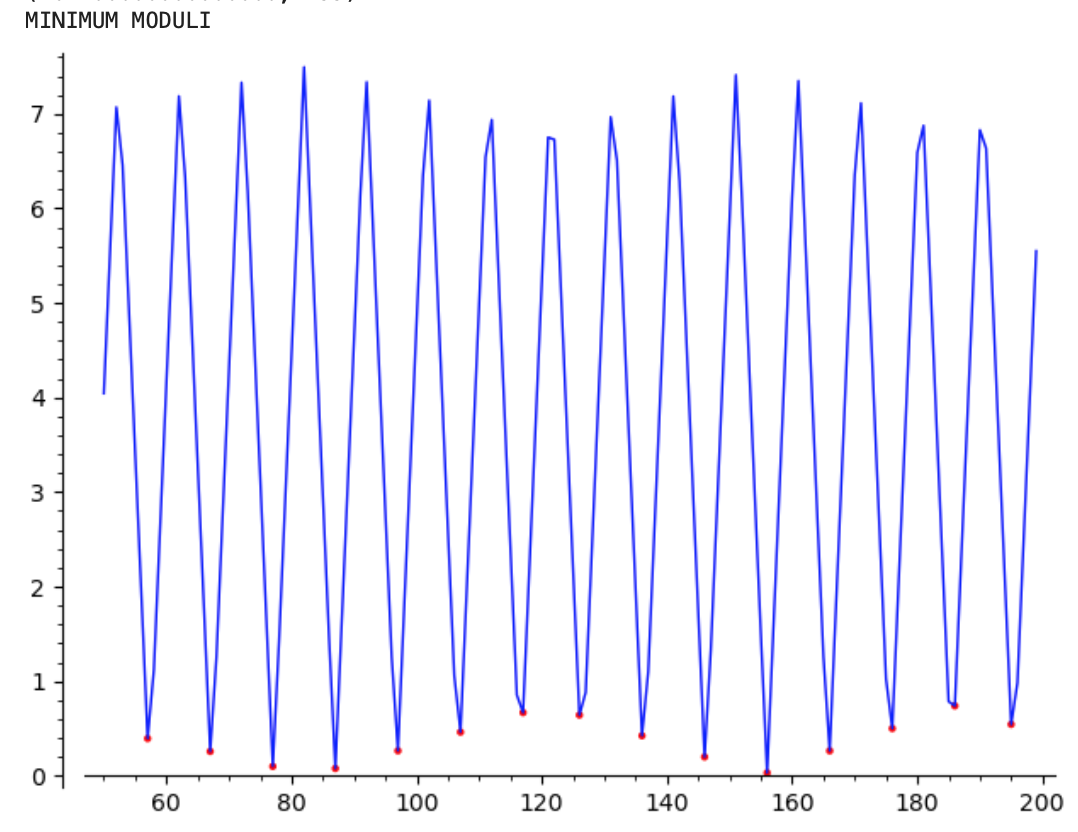}
        \caption*{$v = -.2$}
        \label{fig:v=-.2}
    \end{subfigure}
    
    \begin{subfigure}{0.45\textwidth}
        \centering
        \includegraphics[width=\textwidth]{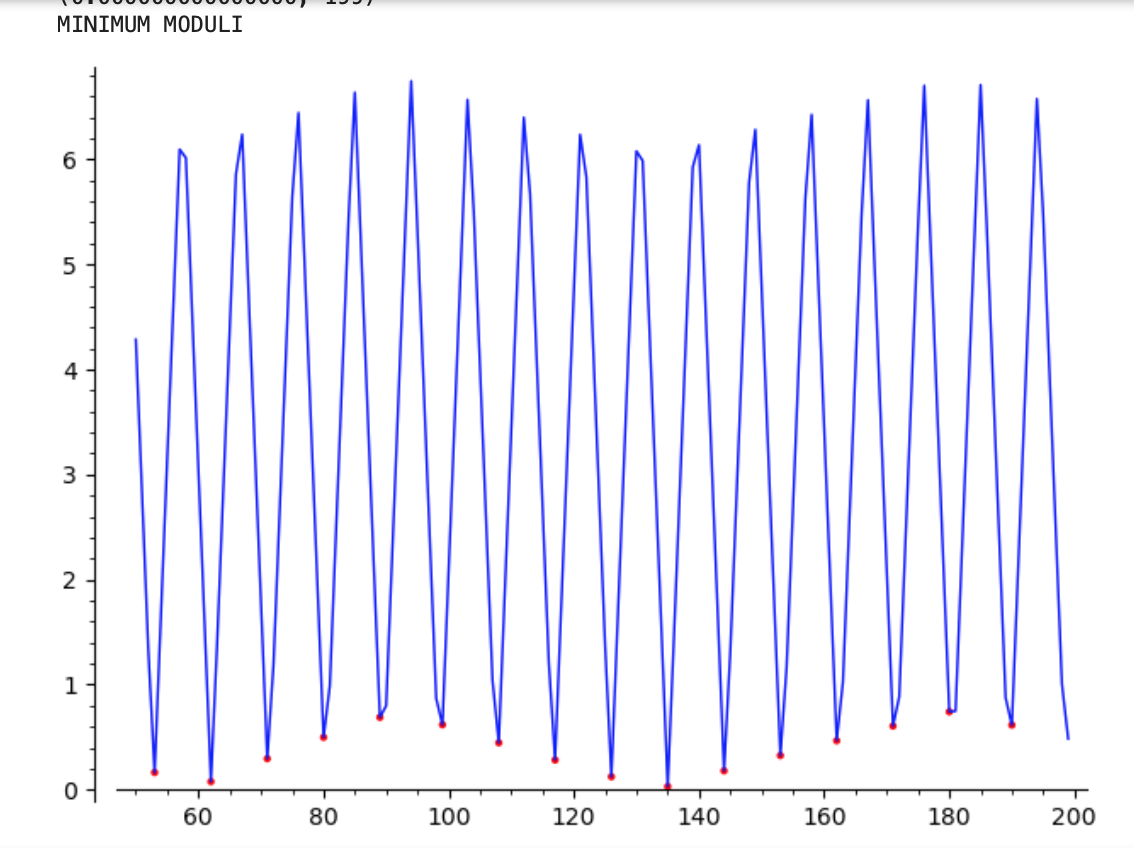}
        \caption*{$v = 0$}
        \label{fig:v=0}
    \end{subfigure}
    \hfill
    \begin{subfigure}{0.45\textwidth}
        \centering
        \includegraphics[width=\textwidth]{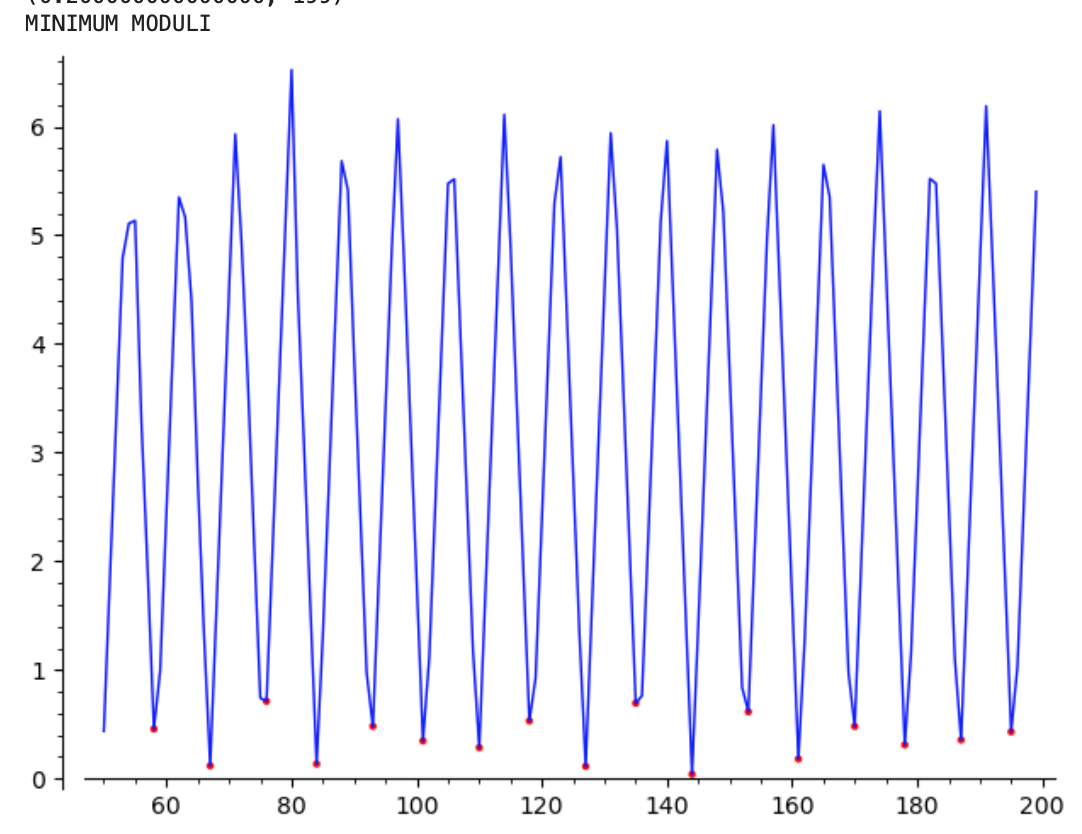}
        \caption*{$v=.2$}
        \label{fig:v=.2}
    \end{subfigure}
     \vspace{1em} % Vertical space between rows
     \begin{subfigure}{0.45\textwidth}
        \centering
        \includegraphics[width=\textwidth]{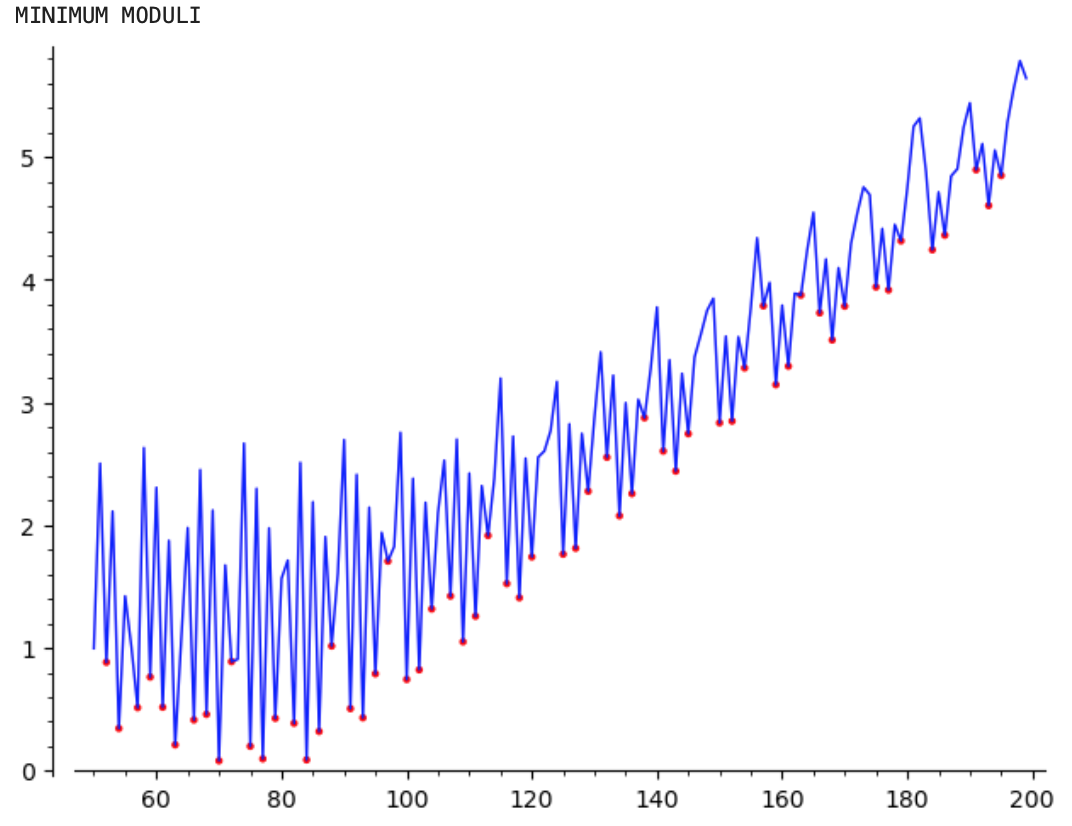}
        \caption*{$v=.4$}
        \label{fig:v=.4}
        \end{subfigure}
     \hfill
    \begin{subfigure}{0.45\textwidth}
        \centering
        \includegraphics[width=\textwidth]{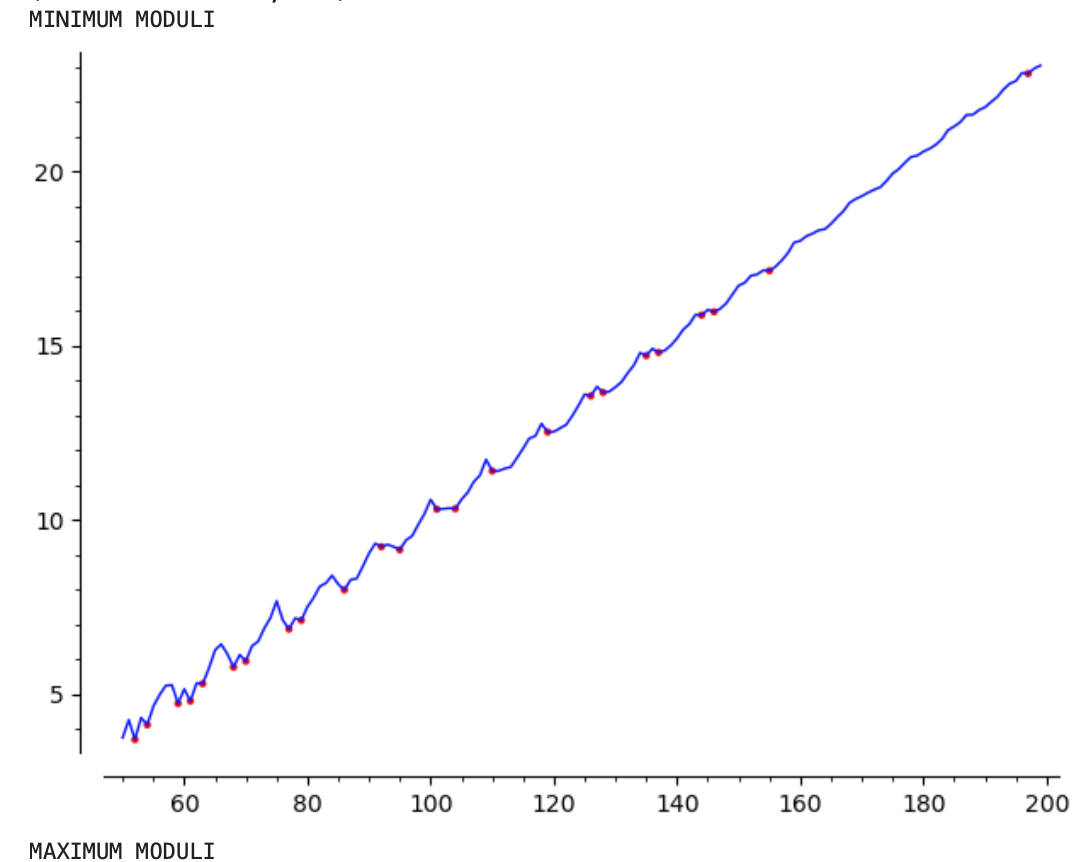}
        \caption*{$v=.6$}
        \label{fig:v=.6}
    \end{subfigure}
    \caption{Minimum moduli for $\mu_v(n)=\mu(n)-v$
    using
    deformed matrices with $c=1$.}
    \label{fig:main}
\end{figure}
\subsection{Color-coded root sets for some integer sequences}\footnote{See the following notebooks in \cite{Bre25}: \texttt{undeformed primeTau 18oct25.ipynb}, \texttt{deformed tauprime 17oct25.ipynb}, \texttt{undeformed tau 18oct25.ipynb}, \path{deformed tau 18oct25.ipynb}, \texttt{deformed $p(n)$ 18oct25.ipynb}, \texttt{undeformed $p(n)$ 19oct25.ipynb}, \texttt{undeformed $p(p_n)$ 18oct25.ipynb}, and \texttt{deformed $p(p_n)$ 18oct25.ipynb}.}
Figure 9 is a gallery of color-coded plots of the roots of the characteristic polynomials associated to the sequences referenced in the captions. The point colors denote the values of $n$ in the relevant sequence to which the given roots belong. Larger $n$ in a sequence will over-write with a new color any root that is too close by. The plots appear to us to indicate that roots tend to lie along certain curves, which would seem to be targets for further investigation.
\begin{figure}[htbp]
    \centering
    \begin{subfigure}{0.45\textwidth}
        \centering
        \includegraphics[width=\textwidth]
        {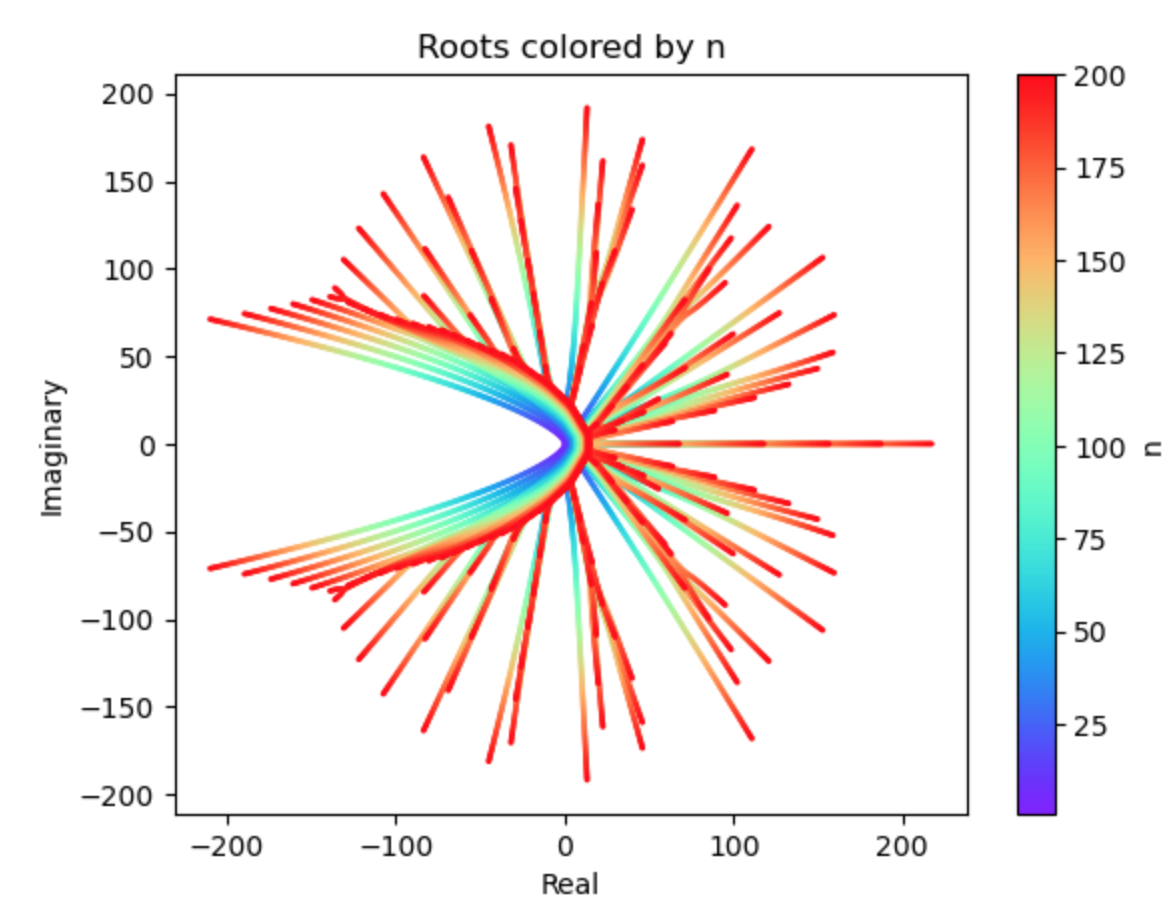}
        \caption*{deformed $p(n)$ roots}
        \label{fig:deformed_p(n)_roots}
    \end{subfigure}
    \hfill
    \begin{subfigure}{0.45\textwidth}
        \centering
        \includegraphics[width=\textwidth]{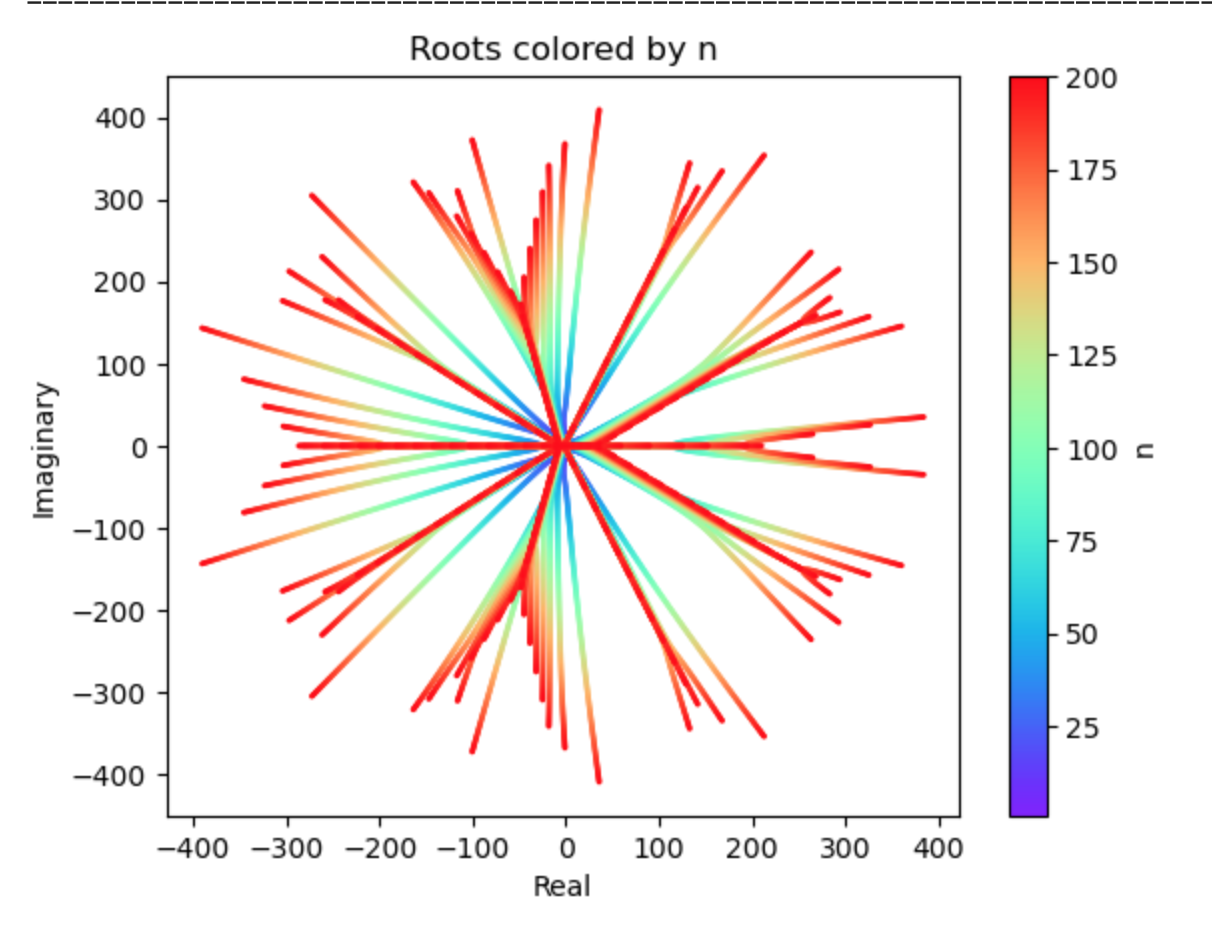}
        \caption*{deformed $p(p_n)$ roots}
        \label{fig:deformed_p_p_n__roots}
    \end{subfigure}
    \begin{subfigure}{0.45\textwidth}
        \centering
        \includegraphics[width=\textwidth]{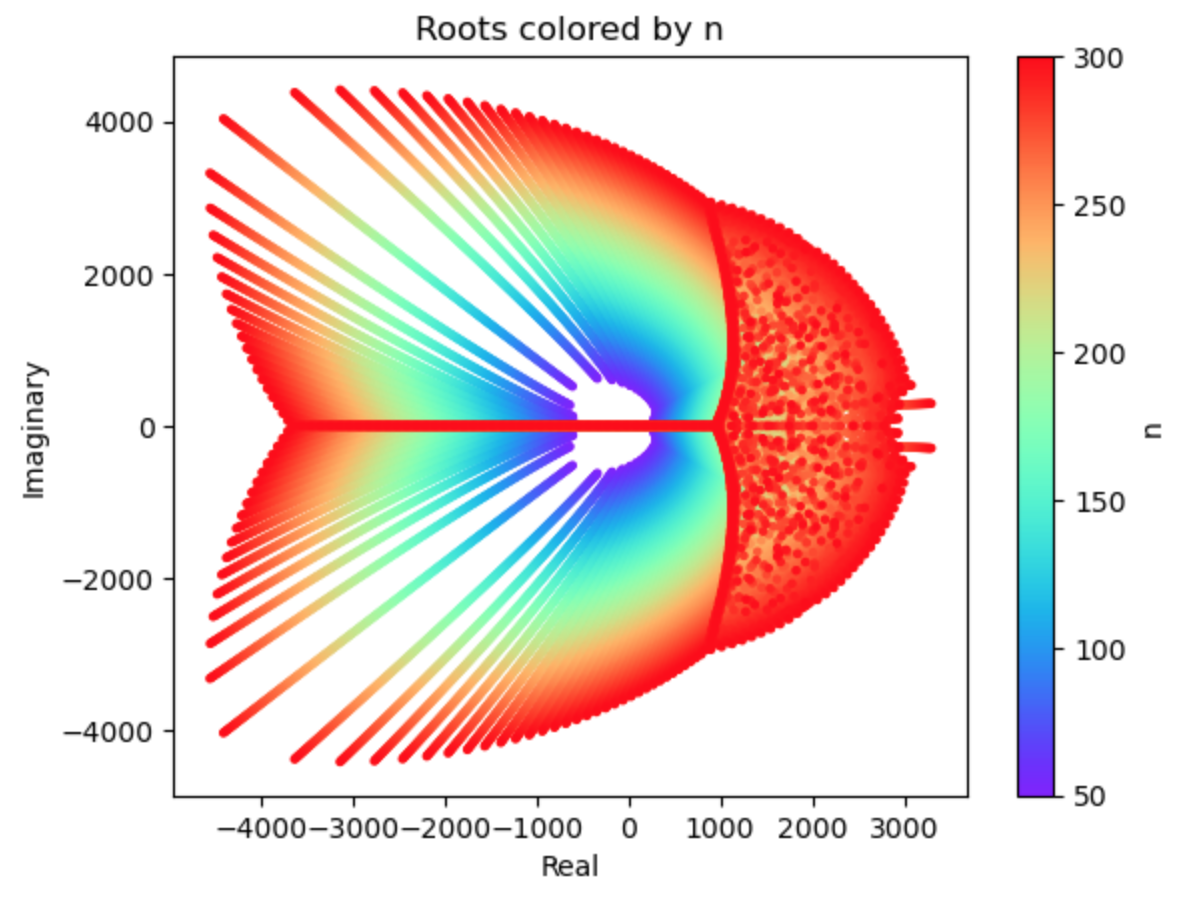}
        \caption*{deformed $\tau(p_n)$ roots}
        \label{fig:deformed_primeTau_roots}
    \end{subfigure}
    \hfill
    \begin{subfigure}{0.45\textwidth}
        \centering
        \includegraphics[width=\textwidth]{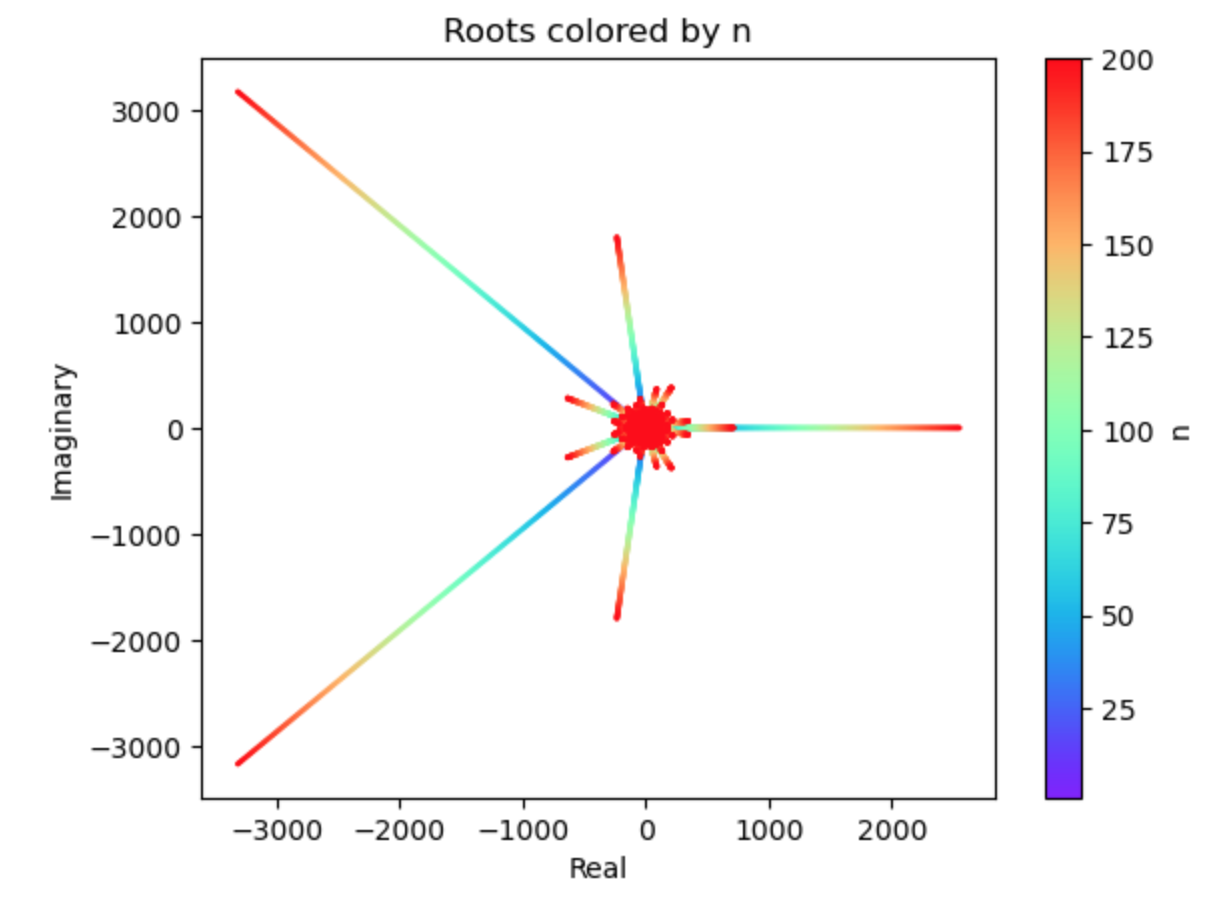}
        \caption*{undeformed $\tau(p_n)$ roots}
        \label{fig:undeformed_primeTau_roots}
    \end{subfigure}
     \vspace{1em} % Vertical space between rows
     \begin{subfigure}{0.45\textwidth}
        \centering
        \includegraphics[width=\textwidth]{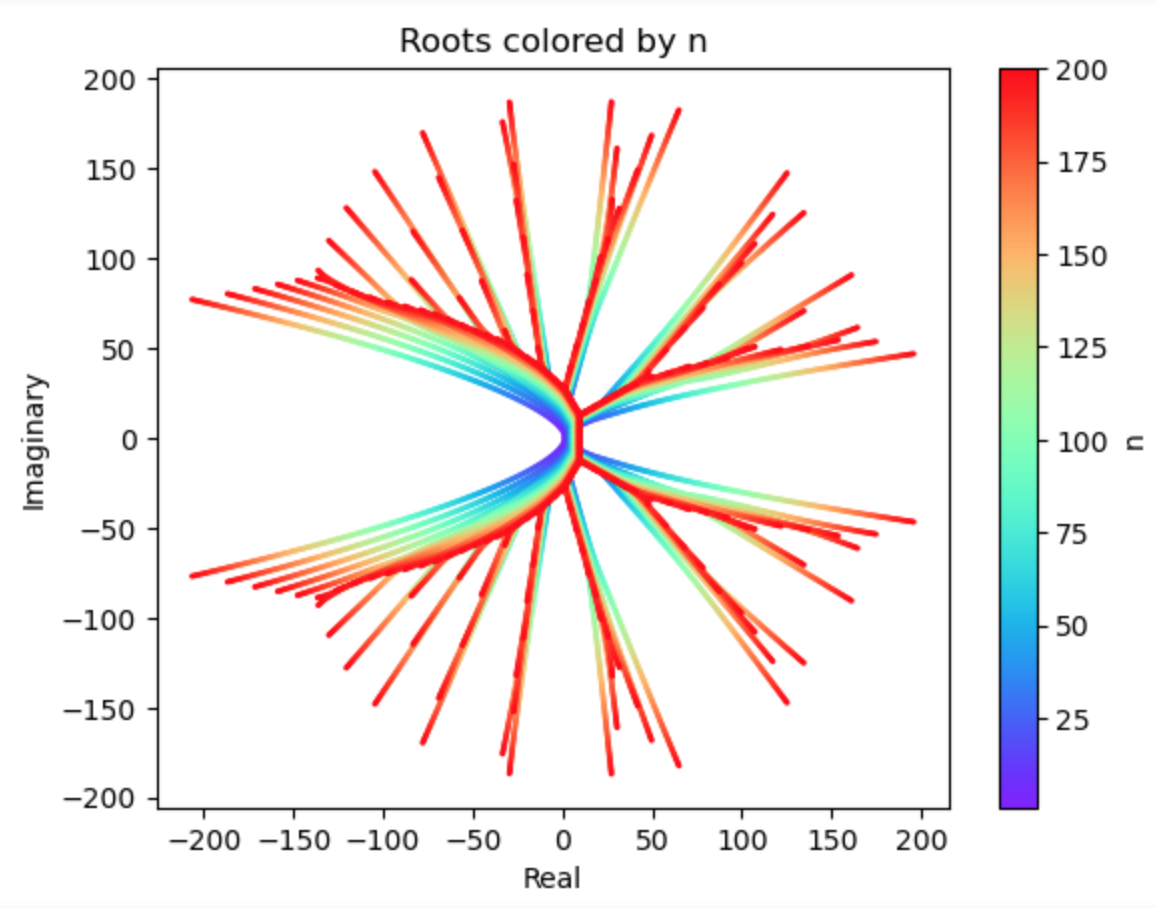}
        \caption*{undeformed $p(n)$ roots}
        \label{fig:undeformed_p_n__roots}
        \end{subfigure}
     \hfill
    \begin{subfigure}{0.45\textwidth}
        \centering
        \includegraphics[width=\textwidth]{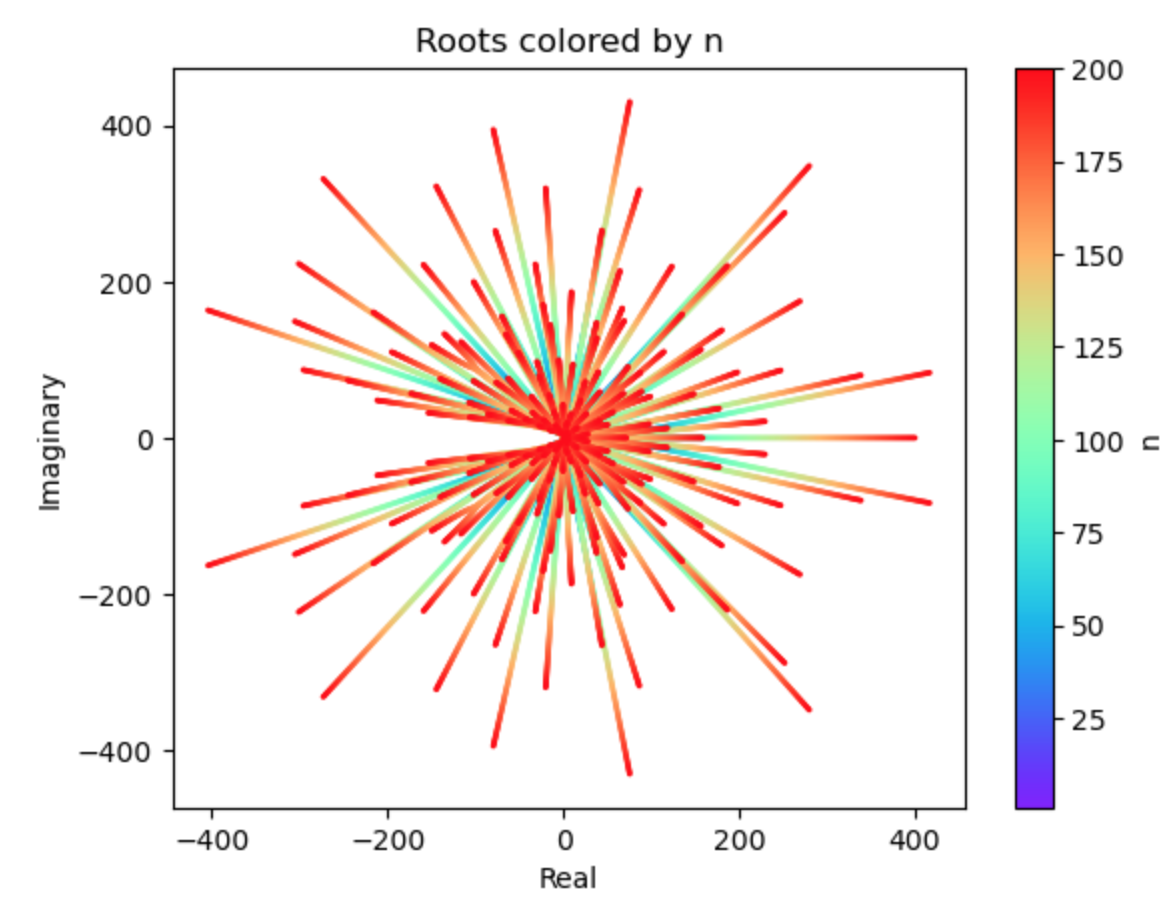}
        \caption*{undeformed $p(p_n)$ roots}
        \label{fig:undeformed_p_p_n}
    \end{subfigure}
     \vspace{1em} % Vertical space between rows
     \begin{subfigure}{0.45\textwidth}
        \centering
        \includegraphics[width=\textwidth]{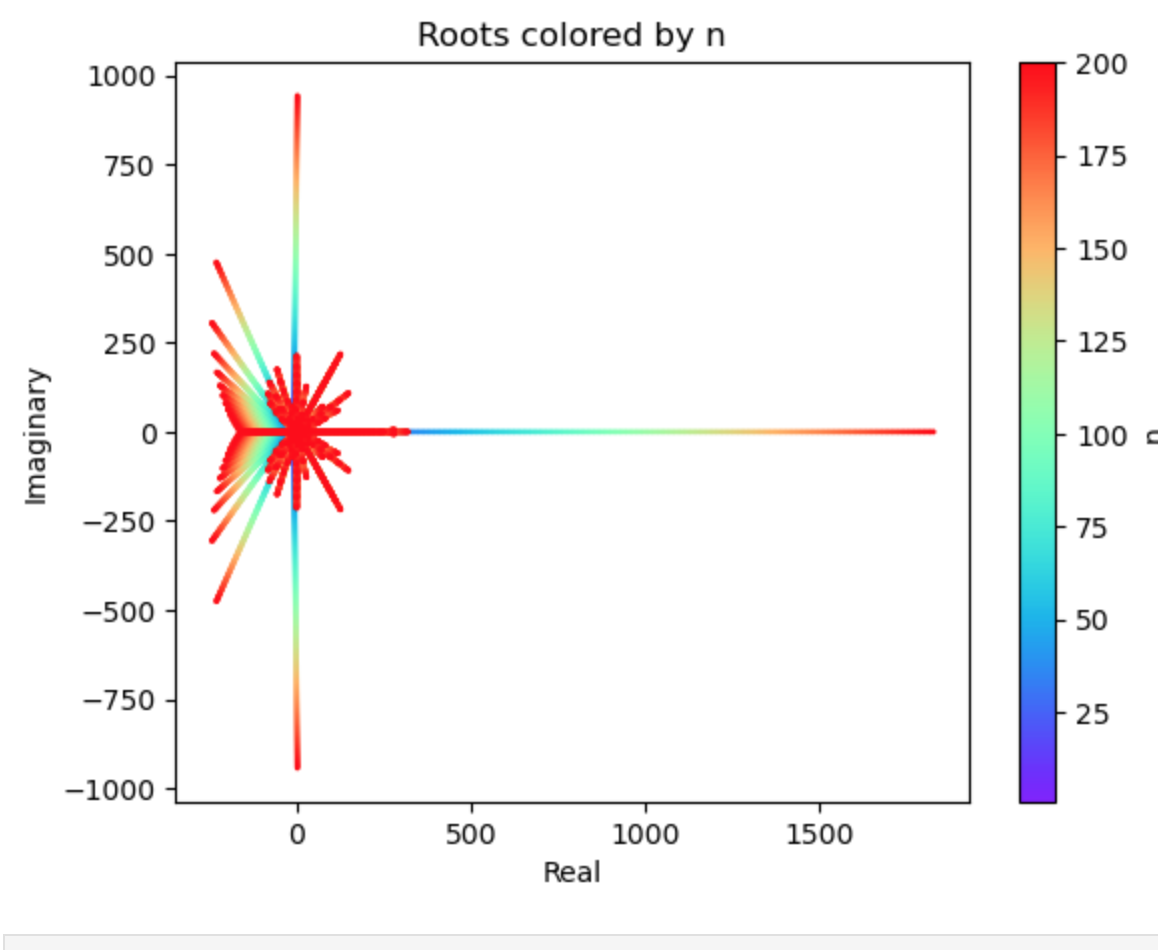}
        \caption*{undeformed tau roots}
        \label{fig:undeformed tau roots}
        \end{subfigure}
     \hfill
    \begin{subfigure}{0.45\textwidth}
        \centering
        \includegraphics[width=\textwidth]{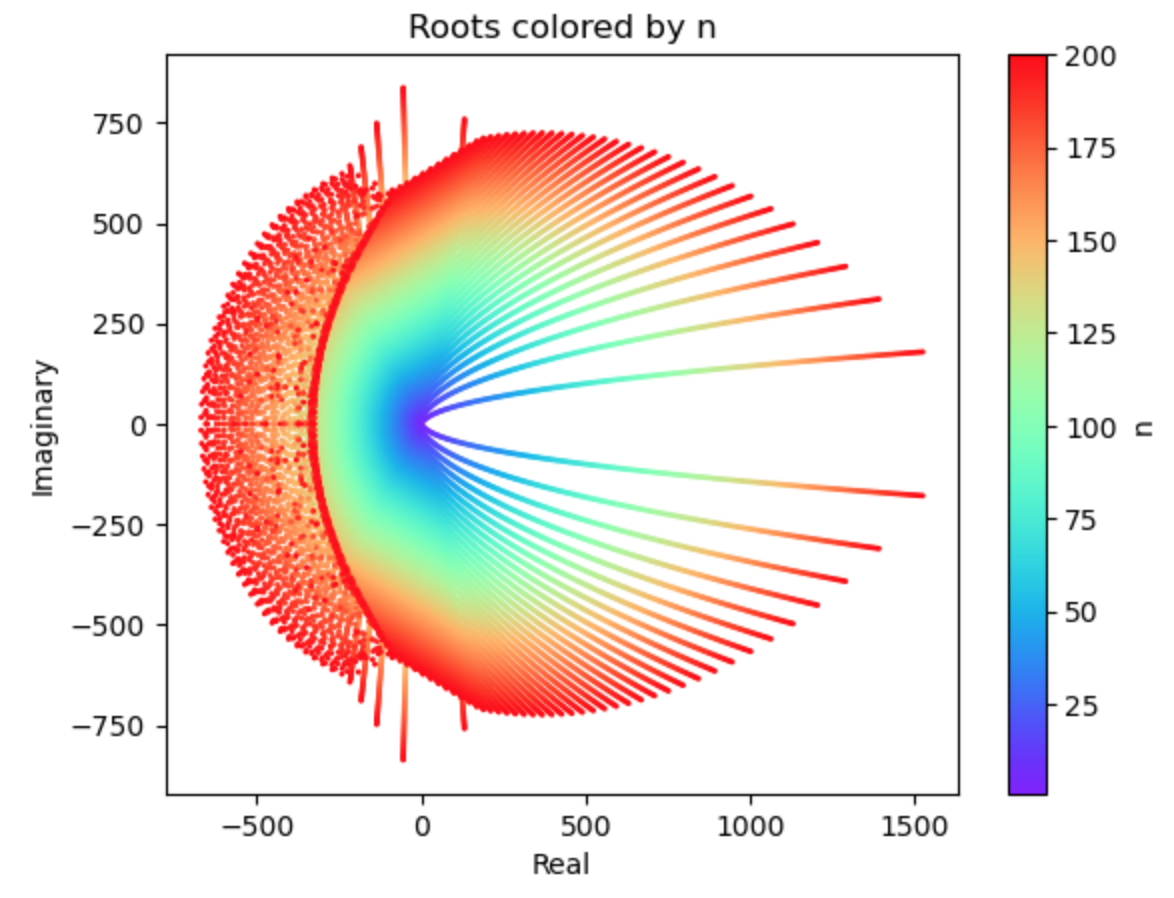}
        \caption*{deformed tau roots}
        \label{fig:deformed_tau_roots}
    \end{subfigure}
    \caption{Ensembles of root sets for characteristic polynomials attached to the $\overline{j}$, color coded by $n$.}
    \label{fig:main2}
\end{figure}

\end{document}